\documentclass[twocolumn]{autart}    % Enable this line and disable the
\usepackage{graphicx}          % Include this line if your
\usepackage{amsmath,amssymb,amsfonts}
\usepackage{enumerate}
\usepackage{color} 
\usepackage{algorithm}
\usepackage{algpseudocode}
\usepackage{cite}
\usepackage{textcomp}
\usepackage{appendix}
\usepackage{diagbox}
\usepackage{placeins}
\usepackage{float}
\usepackage{setspace}
\usepackage{caption}

\newtheorem{mydefn}{Definition}

\newtheorem{myrem}{Remark}

\newtheorem{mylem}{Lemma}
\newtheorem{mytheo}{Theorem}

\usepackage{lipsum}

\makeatletter
\newenvironment{balgorithm}
{% \begin{balgorithm}
	
	\begin{flushleft}
		
		\refstepcounter{algorithm}% New algorithm
		\hrule height.8pt depth0pt \kern2pt% \@fs@pre for \@fs@ruled
		\renewcommand{\caption}[2][\relax]{% Make a new \caption
			{\raggedright\textbf{Pseudocode~\thealgorithm} ##2\par}%
			\ifx\relax##1\relax % #1 is \relax
			\addcontentsline{loa}{algorithm}{\protect\numberline{\thealgorithm}##2}%
			\else % #1 is not \relax
			\addcontentsline{loa}{algorithm}{\protect\numberline{\thealgorithm}##1}%
			\fi
			\kern2pt\hrule\kern2pt
		}
	}{% \end{balgorithm}
		\kern2pt\hrule\relax% \@fs@post for \@fs@ruled
	\end{flushleft}
}
\makeatother

\begin{document}

\begin{frontmatter}

\title{Interval-driven discrete-time general nonlinear robust control: stabilization with closed-loop robust DOA enlargement  \thanksref{footnoteinfo}} % Title, preferably not more

\thanks[footnoteinfo]{This paper was not presented at any IFAC
meeting. Corresponding author Y. Li. Tel. +086-15657191030.}

\author[ZJUT]{Chaolun Lu}\ead{luke@zjut.edu.cn}, 
\author[ZJUT]{Yongqiang Li}\ead{yqli@zjut.edu.cn},
\author[ZJUT]{Zijun Feng}\ead{2111903044@zjut.edu.cn},
\author[QDU]{Zhongsheng Hou}\ead{zshou@qdu.edu.cn}, 
\author[ZJUT]{Yu Feng}\ead{yfeng@zjut.edu.cn}, 
\author[ZJUT]{Yuanjing Feng}\ead{fyjing@zjut.edu.cn},

\address[ZJUT]{College of Information Engineering, Zhejiang University of Technology, Hangzhou, China}
\address[QDU]{School of Automation, Qingdao University, Qingdao, China}

                             % keyword list or with the
                                          % help of the Automatica
                                          % keyword wizard

\begin{abstract}                          % Abstract of not more than 
	This paper presents new results that allow one to address the discrete-time general nonlinear robust control problem. The uncertain system is described by a general nonlinear function set characterized by the nominal model and the corresponding modeling error bound. Traditional synthesis methods design parameters of a structured robust controller. The key aim of this paper is to find an unstructured robust controller set in the state-control space, which enlarges the estimate of the closed-loop robust domain of attraction (RDOA). Based on the interval analysis arithmetic, a numerical method to estimate the unstructured robust controller set is proposed and the rigorous convergence analysis is given. The existing RDOA results are constrained by the level-set of the Lyapunov function, whereas the results in this paper remove this limitation. Furthermore, a solvable optimization problem is formulated so the estimate of RDOA is enlarged by selecting a Lyapunov function from a Lyapunov function set of sum-of-squares polynomials. The method is then validated by a specific case simulation study and results show more extensive RDOA than the previous methods.
\end{abstract}

\end{frontmatter}

 \section{Introduction}
The study of the robust control theory which began in the 1960s is one of the most important branches of the modern control theory. In the early 1970s, in order to deal with some unexpected failures due to the differences between mathematical models and reality, the dominant focus of research shifted from optimality to robustness~\cite{Safonov:2012}. When the available models are no longer "sufficiently accurate", the uncertainty describing how the "true" nominal model might differ from the plant played a big role in the design of the controller. At first, control scientists were generally led to assume that the nominal model is linear~\cite{Petersen:2014,BBhattacharyya:2017}. The main drawback of linear robust control is that the controllers are very conservative when the nonlinearities are significant. That is one of the strong motivations for the development of nonlinear robust control such as the Lyapunov min-max approach~\cite{Corless:1993}, the nonlinear H$_\infty$ approach~\cite{Basar:1995}, the input-to-state stability approach \cite{Sontag:1995_351} and the robust backstepping approach \cite{Freeman:2008}. Published nonlinear robust control approaches almost only consider the nominal models which are affine with respect to the control input. Hence, available information about existing non-affine nonlinearities is ignored. The robust control problem considering the general nonlinear nominal model is a challenge. 

For general nonlinear system without uncertainty, it is very hard to achieve the global stabilization. Hence, the domain of attraction (DOA) of the closed-loop which is an invariant set characterizing asymptotically stabilizable area around the equilibrium point has gained more and more popularity for different systems \cite{Chen:2015_1314,Gering:2015_2231,Li:2014_79,Han:2015}. It is well known that DOA plays an important role both in analysis and synthesis. When the nominal model in robust control is general nonlinear, in order to remain functional despite large changes, we believe that the robust DOA (RDOA) of the closed-loops does require more investigation under certain conditions. For system analysis, considers a given Lyapunov candidate,\cite{Swiatlak:2015} proposes an interval arithmetic approach to obtain the estimation of the RDOA, \cite{Goldsztejn:2019_371} proposes a method of estimating the RDOA for non-smooth systems under uncertainties. For controller design, a few works involve investigating the RDOA of the closed-loop. Only~\cite{Li:2020} proposes a data-driven robust controller design method for general nonlinear discrete systems. The main idea is filtering the simulation data points which can make the Lyapunov function negative definite. However, the method proposed in \cite{Li:2020} has two drawbacks. First, the invariant of the estimate of RDOA is guaranteed by the level-set of the Lyapunov function, which leads to a conservative estimate of RDOA. Second, there is no quantitative analysis result about set estimation errors due to the random sampling and set griding technique approximating sets. A rigorous RDOA estimation method for general systems with uncertainty needs to be developed.

In this paper, considering the discrete-time general nonlinear plant set in which the "true" plant is hidden, the stabilization with closed-loop RDOA enlargement is solved. The plant set is characterized by the general nonlinear nominal model and the corresponding modeling error bound. First, for a given Lyapunov function, the robust negative-definite and invariant set in state-control space (RNIS-SC) is defined. It is shown that the RNIS-SC is an unstructured robust controller set for the plant set, namely, any controller belonging to the RNIS-SC can asymptotically stabilize all the plants in the plant set. Moreover, the projection of the RNIS-SC in the state space is an estimate of the RDOA of the closed-loops. It should be pointed out that the invariant of the RNIS-SC is guaranteed by the infinite iteration of the predecessor set rather than the level-set of the Lyapunov function. However, due to the nonlinearities of the nominal model and the modeling error bound, it is hard to obtain an analytic solution of the RNIS-SC. Then, a numerical method estimating the RNIS-SC for a given Lyapunov function is proposed based on the set estimation via interval analysis (SEVIA) algorithm. The SEVIA algorithm is one of the basic tools in interval analysis arithmetic, which is a kind of numerical method to approximate sets of interest as precise as desired ~\cite{Jaulin:2001}. The convergence of the proposed algorithm approximating RNIS-SC is also proved. Finally, a solvable optimization problem is formulated to enlarge the RDOA of the closed-loops by selecting an appropriate Lyapunov function from a parameterized positive-definite function set.

For the discrete-time general nonlinear robust control problem, the main contributions of this paper consist in: 1) Rather than design parameters of a structured robust controller in traditional synthesis methods, our new synthesis method try to find an unstructured robust controller set (namely, the RNIS-SC) in the state-control space, which enlarges the estimate of closed-loop RDOA and whose invariant is not guaranteed by the level-set of the Lyapunov function. 2) Based on interval analysis arithmetic, a numerical method estimating the RNIS-SC is proposed, for which the rigorous convergence analysis is given.

The rest of this paper is organized as follows. In Section 2, a general nonlinear plant set is formulated and our control objective is proposed. Then, a sufficient condition for asymptotic stabilization and estimation of the RDOA is given with the definition of RNIS-SC. In Section 3, we briefly introduced interval arithmetic, then our main algorithm is proposed. Next, based on the theoretical result in the preceding section, the stabilization problem is solved by using the interval analysis approach to estimate the RNIS-SC. In Section 4, the closed-loop RDOA enlargement method and the controller design method is proposed. In Section 5, the controller design method is verified with the simulation result in a specific case. The conclusion is drawn in Section 6.

\textbf{Notation: } For a vector $x \in \mathbb{R}^n$, $x_{(i)}$ represents the $i$-th element of $x$, $i = 1,2,\cdots,n$. For $x_1,x_2\in\mathbb{R}^n,x_1\leq x_2$ means $x_1$ is less than or equal to $x_2$ element by element. For two vectors $x \in \mathbb{R}^n$ and $u \in \mathbb{R}^m$, $w = (x,u)$ represents a new vector in $\mathbb{R}^{n+m}$. $[w] \subset \mathbb{R}^{n+m}$ represents a box belonging to $\mathbb{R}^{n+m}$ (\textit{e.g., }a rectangular region when $n+m = 2$). $\mathbb{W} \subset \mathbb{R}^{n+m}$ represents an arbitrary compact~\cite{Kreyszig:1978} subset of $\mathbb{R}^{n+m}$. $\hat{\mathbb{W}}$ represents an approximation of $\mathbb{W}$ by covering $\mathbb{W}$ with non-overlapping boxes in a set of boxes (there is no ambiguity when $\hat{\mathbb{W}}$ is viewed as a set of boxes or the unions of boxes in a set of boxes according to contents). sub($x$) represents supremacy of $x$, inf($x$) represents infimum of $x$. $[x] \subset \mathbb{R}^n$ reprensents an interval. $[x]\sqcup[y] = [\min\{\underline{x},\underline{y}\}, \max\{\bar{x},\bar{y}\}]$ represents the interval union of the interval $[x]$ and $[y]$.

\section{Robust Stabilization Controller Set: Robust Negative-definite and Invariant Set in State-Control Space}

Consider the nonlinear discrete-time plant set
\begin{eqnarray}
& \!\! \mathfrak{F} := \Big\{ f: \mathbb{R}^n \!\times\! \mathbb{R}^m \!\to \!	\mathbb{R}^n \Big| f(0,0) = 0, \nonumber \\ 
& \quad \hat{f}(x,u) - \delta(x,u) \leq f(x,u) \leq \hat{f}(x,u) + \delta(x,u) \Big\}, \label{eq:plant_set}
\end{eqnarray}
where $x(k) \in \mathbb{R}^n$ is the state, $u(k) \in \mathbb{R}^m$ is the control input, continuous function $\hat{f}: \mathbb{R}^n \times \mathbb{R}^m \to \mathbb{R}^n$ is the known nominal model satisfying $\hat{f}(0,0) = 0$ and continuous function $\delta: \mathbb{R}^n \times \mathbb{R}^m \to \mathbb{R}^n_+$ is the known modeling error bound satisfying $\delta(0,0) = 0$. The control objective is to find a robust controller $\mu: \mathbb{R}^n \to \mathbb{R}^m$ and an estimate of the robust closed-loop DOA to ensure that, $\forall f \in \mathfrak{F}$, the closed-loop $x(k+1) = f(x(k),\mu(x(k))), k=0,1,2,...,$ is asymptotically stable at the origin for all initial states in the estimate of the closed-loop RDOA and to enlarge the estimate of the RDOA as large as possible.

The main objective of this section consists in characterizing robust negative-definite and invariant set in the state-control space (RNIS-SC), which is an unstructured robust stabilization controller set, and discussing properties of RNIS-SC. Moreover, it is observed how such features can be employed to obtain the estimation of RDOA.

\subsection{Robust Negative-definite Set in State-Control Space}

For a given nominal model $\hat f$ and the known modeling error bound $\delta$, we omit time instant $k$ and let $x_+, x$ and $u$ denote $x(k+1), x(k)$ and $u(k)$, respectively. Then the dynamic of plant set \eqref{eq:plant_set} can be represented as a domain $\Pi$ in $(2n+m)$-dimensional space, which is defined as 
\begin{eqnarray}
\Pi= {\rm{\{ (}}{x_ + },x,u{\rm{)}} \in \mathbb{R}^{2n+m} |\hat f(w) - \delta (w) \le {x_ + } \nonumber \\ \le \hat f(w) + \delta (w),
 w = {(x,u)}{\rm{\} }}.
\end{eqnarray}
For a given Lyapunov function $L:\mathbb{R}^{n}\rightarrow\mathbb{R}$, the subset $\tilde{\Pi}(L)$ of $\Pi$ is defined as
\begin{eqnarray}
\tilde \Pi_N(L) = \{ {\rm{(}}{x_ + },x,u{\rm{)}} \in {\Pi_{\mathfrak{F}}}|L({x_ + }) - L(x) < 0\}.
\end{eqnarray} 
All points $(x_+,x,u)$ belong to the $\tilde \Pi_N(L)$ make the time difference of the Lyapunov function $L$ negative definite, so $\tilde \Pi_N(L)$ is called negative-definite set in state-state-control space (ND-SSC). However, ND-SSC $\tilde \Pi_N(L)$ is not robust for plant set~\eqref{eq:plant_set}. The concept of the robust negative-definite set is proposed in \cite{Li:2020}. Any point $(x_+,x,u)$ in the robust negative-definite set must satisfies: $\forall f\in\mathfrak{F}, L(f(x,u))-L(x)<0$. Hence, robust negative-definite set in state-state-control space (RNS-SSC) $\Pi_N(L)$ is
\begin{eqnarray}
\Pi_N(L)=&\{(x_+,x,u)\in \tilde \Pi_N(L)|\forall x_+' \in \mathbb{X}_+(x,u),\nonumber
	\\&L(x_+')-L(x)<0,x_+\in\mathbb{X}_+(x,u)\},
\end{eqnarray}
where, for a given $(x,u)\in \mathbb{R}^{n+m}$, $\mathbb{X}_+(x,u)$ is defined as
\begin{eqnarray}
\mathbb{X}_+(x,u)=& \{x_+\in \mathbb{R}^{n}|\hat{f}(w)-\delta(w)\leq x_+ \leq \nonumber\\
&\hat{f}(w)+\delta (w), w=(x,u)\}. \label{eq:nextset_x}
\end{eqnarray}
The robust negative-definite set in state-control space can be obtained by projecting RND-SSC along the $x_+$ space to the state-control space which is defined as follows.

\begin{mydefn}
	For a given Lyapunov function $L$, the robust negative-definite set $\mathbb{W}_{\rm N}(L)$ for the plant set~\eqref{eq:plant_set} is defined as
	\begin{eqnarray}
	\mathbb{W}_{\rm N}(L)= & \{(x,u) \in \mathbb{R}^{n+m}|\forall x_+ \in\mathbb{X}_+(x,u),\nonumber\\
	& L(x_+)-L(x)<0 \}, \label{def:w_nl}
	\end{eqnarray}
	where $\mathbb{X}_+(x,u)$ is defined in~\eqref{eq:nextset_x}. \label{def:negative definite}
\end{mydefn}

\begin{mylem}\label{Lem:RNS}
	For plant set~\eqref{eq:plant_set}, if a Lyapunov function $L$ and a controller $\mu$ exist such that
	\begin{equation}
	\forall x \in {\rm{proj}}({\mathbb{W}_{\rm{N}}}(L)),(x,\mu (x)) \in {\mathbb{W}_{\rm{N}}}(L),0 = \mu (0)\nonumber
	\end{equation}
	then $L\left( {x(k+1)} \right) - L\left( {x(k)} \right) < 0$ for the closed-loop $x(k+1) = f\left( {x(k),\mu (x(k))} \right),$ when $x(k)\in {\rm{proj}}({\mathbb{W}_{\rm{N}}}(L)),\forall f\in\mathfrak{F}$.
\end{mylem}

\begin{pf}
	We have $(x(k),\mu (x(k))) \in {\mathbb{W}_{\rm{N}}}(L)$, according to the definition of RNS-SC $\mathbb{W}_{\rm N}(L)$ in~\eqref{def:w_nl}, then
	\begin{equation}
	\forall x(k+1) \in\mathbb{X}_+(x(k),\mu (x(k))),L(x(k+1))-L(x(k))<0. \nonumber
	\end{equation}
	Based on \eqref{eq:nextset_x}, $\forall f\in\mathfrak{F}, f(x(k),\mu (x(k)))\in\mathbb{X}_+$. Hence, for any state ${x(k)} \in {\rm{proj}}({\mathbb{W}_{\rm{N}}}(L))$ and $\forall f \in\mathfrak{F} $, $L(x(k+1))-L(x(k))<0$ is satisfied. \hfill $\blacksquare$
\end{pf}

Based on the Lemma~\ref{Lem:RNS}, one may think that ${\rm proj}(\mathbb{W}_{\rm N}(L))$ can be regarded as an estimate of closed-loop RDOA. However, it is wrong because there is no guarantee that the next state $x(k+1)$ is still in ${\rm proj}(\mathbb{W}_{\rm N}(L))$. Once the next state $x(k+1)$ does not belong to ${\rm proj}(\mathbb{W}_{\rm N}(L))$, the time difference of Lyapunov function  is no longer negative definite. Finding a level set to deliver an enclosure of the Lyapunov function is one of the solutions~\cite{Li:2014_79,Swiatlak:2015,Goldsztejn:2019_371} but leads to the conservative result.
 
\subsection{Robust invariant set in state-control space} 
 
For the purpose of finding an unstructured robust controller set, we give the definition of the robust invariant set in the state-control space (RIS-SC) as follows.

\begin{mydefn}
	The robust invariant set $\mathbb{W}_{\rm I}$ for the plant set~\eqref{eq:plant_set} is defined as
	\begin{eqnarray}
	\mathbb{W}_{\rm I}=\{w=(x,u)\in\mathbb{R}^{n+m}|\mathbb{X}_+(w)\subset{\rm proj}(\mathbb{W}_{\rm I})\}, \label{eq:W_I}
	\end{eqnarray} 
	where $\mathbb{X}_+$ is defined in~\eqref{eq:nextset_x}, ${\rm proj}(\mathbb{W}_{\rm I})\subset\mathbb{R}^{n}$ represents the prejection of $\mathbb{W}_{\rm I}\subset\mathbb{R}^{n+m}$ along the state space $\mathbb{R}^n$. \label{def:invariant}
\end{mydefn}

\begin{mylem} \label{Lem:RIS}
	For the plant set~\eqref{eq:plant_set}, if a controller $\mu$ satisfies
	\begin{equation}
	\forall x \in {\rm{proj}}({\mathbb{W}_{\rm{I}}}),(x;\mu (x)) \in {\mathbb{W}_{\rm{I}}}, \label{eq:xinW_I}
	\end{equation}
	then the solution $\phi ({x_0},k)$ of the closed-loop $x(k + 1) = f\left( {x(k),\mu (x(k))} \right)$ with any initial state ${x_0} \in {\rm{proj}}({\mathbb{W}_{\rm{I}}})$  is in ${\rm{proj}}({\mathbb{W}_{\rm{I}}})$ for all $k>0$.
\end{mylem}

\begin{pf}
	We prove Lemma~\ref{Lem:RIS} in the form of mathematical induction.
	
	For $k=0$, we have $\phi(x_0,0)=x_0 \in {\rm proj}(\mathbb{W}_{\rm I})$. Next, it is proofed that $\phi(x_0,k+1) \in {\rm proj}(\mathbb{W}_{\rm I})$ when $\phi(x_0,k) \in {\rm proj}(\mathbb{W}_{\rm I})$. Because $x(k)=\phi(x_0,k) \in {\rm proj}(\mathbb{W}_{\rm I})$, from~\eqref{eq:xinW_I}, $(x(k);\mu(x(k))) \in \mathbb{W}_{\rm I}$ is obtained. According to Definition~\ref{def:invariant}, $ \forall f \in\mathfrak{F}$, we have $x(k+1)=\phi(x_0,k+1) \in {\rm proj}(\mathbb{W}_{\rm I}).$
	\hfill $\blacksquare$
\end{pf}

 With the definition of RIS-SC, robust negative-definite and invariant set in state-control space are introduced in the next subsection. 

\subsection{Robust negative-definite and invariant set in state-control space}

For a given Lyapunov function $L$, combine the definition of the robust negative-definite set $\mathbb{W}_{\rm N}(L)$ and the robust invariant set $\mathbb{W}_{\rm I}$, the definition of the robust negative-definite and invariant set in the state-control (RNIS-SC) space are introduced as follows.

\begin{mydefn} \label{defn:NDIset}
	The robust negative-definite and invariant set $\mathbb{W}_{\rm N\&I}(L)$ for the plant set~\eqref{eq:plant_set} is defined as 
\begin{eqnarray}
  \mathbb{W}_{\rm N\&I}(L)=&\{w=(x,u)\in\mathbb{R}^{n+m}|\nonumber
  \\&\mathbb{X}_+(w)\subseteq {\rm proj}(\mathbb{W}_{\rm N\&I}(L))\}, \label{eq:ri}\\
  &\forall x_+\in \mathbb{X}_+, {L(x_+)}-L(x)<0, \label{eq:rnd}
\end{eqnarray}
where $x_+$ is the state at the next time instant, $\mathbb{X}_+$ is defined in~\eqref{eq:nextset_x}, $L(x)$ is a given Lyapunov function. 
\end{mydefn}

Since $\mathbb{W}_{\rm N\&I}(L)$ is robust negative definite, according to the~\eqref{def:w_nl}, $\forall (x(k),u(k))\in\mathbb{W}_{\rm N\&I}(L), \forall f\in\mathfrak{F}$, the time difference of the Lyapunov function $L(x(k))$ at $k$ instant is negative definite. Futher more, $\mathbb{W}_{\rm N\&I}(L)$ is robust invariant, according to the~\eqref{eq:W_I}, $\forall (x(k),u(k))\in\mathbb{W}_{\rm N\&I}(L),\forall f\in\mathfrak{F}$, the next state $x(k+1)=f(x(k),u(k))$ is still in the ${\rm proj}(\mathbb{W}_{\rm N\&I}(L))$, which means that there exists $u(k+1)$ such that the time difference of the Lyapunov $L(x(k+1))$ is still negative difinite at $k+1$. Hence, the closed-loop system is asymptotically stable for all initial state in ${\rm proj}(\mathbb{W}_{\rm N\&I}(L))$ which is shown in the following theorem.

\begin{mytheo} \label{Theo:WNI}
	Given RNIS-SC $\mathbb{W}_{\rm N\&I}(L)\subset\mathbb{R}^{n+m}$ for the plant set~\eqref{eq:plant_set} and Lyapunov function $L:\mathbb{R}^n\rightarrow\mathbb{R}_+$, if the controller $\mu:\mathbb{R}^n\rightarrow\mathbb{R}^m$ satisfies
	\begin{eqnarray}
      \forall x\in{\rm{proj}}({\mathbb{W}_{\rm{N\&I}}(L)}),(x;\mu(x))\in\mathbb{W}_{\rm N\&I}(L),\nonumber\\0=\mu(0), \label{eq:proposi_suff}
	\end{eqnarray}
	then the closed-loop system $x(k+1)=f(x(k),\mu(x(k))), $ $\forall f \in \mathfrak{F}$ is asymptotically stable for any initial state $x_0 \in {\rm{proj}}({\mathbb{W}_{\rm{N\&I}}})$.
\end{mytheo}
\begin{pf}
	For any $f\in\mathfrak{F}$, let $\phi(x_0,k)$ denote the solution of $x(k+1) = f\big(x(k),\mu(x(k))\big)$ at time $k$ with the initial state $x_0$. According to the Definition~\ref{defn:NDIset}, $\mathbb{W}_{\rm N\&I}(L)$ is robust negative-definite. From~\eqref{eq:proposi_suff} and Lemma~\ref{Lem:RNS}, we have, $\forall f\in\mathfrak{F}$,	\begin{equation}
	 L(\phi(x_0,k+1)) - L(\phi(x_0,k)) <0, \label{eq:prop:stab:pf:L(k+1)<L(k)}
	\end{equation}
	when $\phi(x_0,k)\in{\rm proj}(\mathbb{W}_{\rm N\&I}(L))$. According to the Definition~\ref{defn:NDIset}, $\mathbb{W}_{\rm N\&I}(L)$ is also robust invariant. From~\eqref{eq:proposi_suff} and Lemma~\ref{Lem:RNS}, we have, $\forall f\in\mathfrak{F}$, $\forall x_0\in {\rm proj}(\mathbb{W}_{\rm N\&I}(L))$, $\forall k\ge 0$,
	\begin{equation}
	\phi(x_0,k)\in {\rm proj}(\mathbb{W}_{\rm N\&I}(L)). \label{eq:eq:prop:stab:pf:phi(x,k)_in_WNI}
	\end{equation}
	From~\eqref{eq:prop:stab:pf:L(k+1)<L(k)} and~\eqref{eq:eq:prop:stab:pf:phi(x,k)_in_WNI}, it follows that, $\forall f\in\mathfrak{F}$, $\forall x_0 \in {\rm proj}(\mathbb{W}_{\rm N\&I}(L))$, $L(\phi(x_0,k))$ is monotonically decreasing with time. And because $L$ is positive-definite, $L(\phi(x_0,k))$ is bounded from below by zero. Hence, 
	\begin{equation}
	\forall x_0 \in {\rm proj}(\mathbb{W}_{\rm N\&I}(L)), \lim_{k \to \infty} L(\phi(x_0,k)) = 0.  \label{eq:reductio}
	\end{equation}
	From the above equation, we can derive that
	\begin{equation}
	\forall x_0 \in {\rm proj}(\mathbb{W}_{\rm N\&I}(L)), \lim_{k \to \infty} \phi(x_0,k) = 0, \label{eq:absurdum}
	\end{equation}
	The process of \eqref{eq:reductio} to \eqref{eq:absurdum} is omitted. This can be proven by reductio ad absurdum (details see
	the proof of Theorem 13.2 in \cite{Haddad:2008}).
	\hfill $\blacksquare$
\end{pf}

With the Theorem~\ref{Theo:WNI}, once the RNIS-SC $\mathbb{W}_{\rm N\&I}(L)$ is obtained, any controller satisfies~\eqref{eq:proposi_suff} can stabilize the plant set~\eqref{eq:plant_set}. However, due to nonlinearities of $\hat{f}$ and $\delta$, it is hard to find the analytic solution of $\mathbb{W}_{\rm N\&I}(L)$. In Section 4, we propose a numerical method to obtain an approximate solution of $\mathbb{W}_{\rm N\&I}(L)$. It is obvious that RNIS-SC $\mathbb{W}_{\rm N\&I}(L)$ is an invariant subset of RNS-SC $\mathbb{W}_{\rm N}(L)$. In the next subsection, a theoretical result about how to find an invariant subset of RNS-SC is given, which is the theoretical cornerstone of the numerical method proposed in Section 4.

\subsection{Finding RNIS-SC by predecessor operator}

Because RNIS-SC $\mathbb{W}_{\rm N\&I}(L)$ is an invariant subset of RNS-SC $\mathbb{W}_{\rm N}(L)$, we first give the theoretical result about finding an invariant subset of any interested region in the state-control space. Predecessor set (or one-step backward set) of a given set is the center of computing invariant sets and is well studied in~\cite{Rakovic:2006,Bertsekas:1972,Blanchini:1999}. While the predecessor set in existing literature is defined in the state space, we define the predecessor set in the state-control space as follows.

	\begin{equation}
	{\rm{Pre}}(\mathbb{W}) = \left\{ {w \in \mathbb{R}{^{n + m}}\left| {{\mathbb{X}_ + }(w) \subset {\rm{proj}}(\mathbb{W})} \right.} \right\} \label{eq:defn:Predecessor}
	\end{equation}
	
With predecessor set ${\rm{Pre}}(\mathbb{W})$ for $\mathbb{W}\subset \mathbb{R}{^{n + m}}$ , we define a new mapping between subsets of the state-control space as follows.
\begin{mydefn}
	\label{def:map}
	For $\mathbb{W} \subset \mathbb{R}^{n+m}$, mapping $\mathcal{I}$ is defined as
	\begin{equation}
	\mathcal{I}(\mathbb{W}) ={\rm{Pre}}(\mathbb{W})\cap\mathbb{W} = \Big\{w \in \mathbb{W}\Big|\mathbb{X}_+(w) \subset \mathrm{proj}(\mathbb{W})\Big\}, \label{eq:defn:I}
	\end{equation}
	where $\mathbb{X}_+$ is shown in~\eqref{eq:nextset_x}, $\mathrm{proj}(\mathbb{W})\subset\mathbb{R}^n$ represents the orthogonal projection of $\mathbb{W}\subset \mathbb{R}^{n+m}$ along the control space to the state space.
\end{mydefn}

\begin{myrem}
	Mapping $\mathcal{I}$ is similar to the mapping $\mathcal{C}$ defined in (11) of \cite{Bertsekas:1972}.
	From \eqref{eq:defn:I}, it is obvious that the fundamental principles of mappings $\mathcal{I}$ and $\mathcal{C}$ are similar. The difference is that $\mathcal{I}$ maps subsets of the state-control space into subsets of the state-control space, while $\mathcal{C}$ maps subsets of the state space into subsets of the state-control space.
\end{myrem}

The composition of mapping $\mathcal{I}$ with itself $i$ times is denoted by $\mathcal{I}^i$. Mapping $\mathcal{I}$ has the following property.
\begin{mylem} \label{lem:I_infty}
	For any compact set $\mathbb{W} \subset \mathbb{R}^{n+m}$, set limit $\mathcal{I}^{\infty}(\mathbb{W}) := \lim_{i \to \infty} \mathcal{I}^{i}(\mathbb{W}) := \bigcap_{i = 1}^{\infty} \mathcal{I}^i(\mathbb{W})$ exists and is invariant for plant set~\eqref{eq:plant_set}.
\end{mylem}
\begin{pf}
	From the definition of mapping $\mathcal{I}$, we have
	\begin{eqnarray}
	\mathcal{I}^{i+1}(\mathbb{W}) = \Big\{w \in \mathcal{I}^{i}(\mathbb{W})\Big|f(w) \in \mathrm{proj}(\mathcal{I}^{i}(\mathbb{W}))\Big\} \label{eq:prop:I_infty:pf:I_i}
	\end{eqnarray}
	According to \eqref{eq:prop:I_infty:pf:I_i}, it is obvious that  $\{\mathcal{I}^i (\mathbb{W})\}$ is a monotonically decreasing sequence of sets in the sense that $\mathcal{I}^{i+1} (\mathbb{W}) \subset \mathcal{I}^i (\mathbb{W})$ for all $i \geq 1$. Moreover, $\mathcal{I}^i (\mathbb{W})$ is a closed set for all $i \geq 1$. Hence, the set limit of $\mathcal{I}^i (\mathbb{W})$ exists (see \cite{Rockafellar:2009}, pp. 111).
	
	$\mathbb{W}$ is compact set (which is equivalent to bouneded and closed in metric space \cite{Kreyszig:1978} (p.77, Lemma 2.5-2)). Let $\mathbb{U}$ represents the orthogonal projection of $\mathbb{W}\subset \mathbb{R}^{n+m}$ along the state space to the control space, ${\rm proj}(\mathbb{W})$ and $\mathbb{U}$ are both bounded and $\forall f\in\mathfrak{F}$, $f$  is continuous over $\mathbb{R}^n\times \mathbb{R}^m$, which indicates that $\mathcal{I}(\mathbb{W})$ is bounded. Now we show that  $\mathcal{I}(\mathbb{W})$ is closed. Let ${w_i}={(x_i,u_i)}$ be a convergent sequence in  $\mathcal{I}(\mathbb{W})$. Since ${\rm proj}(\mathbb{W})$ and $\mathbb{U}$ are compact, it follows that $\bar x = {\lim _{i \to \infty }}{x_i} \in \mathbb{X}$ and $\bar u = {\lim _{i \to \infty }}{u_i} \in \mathbb{U}$. Let $y_i =f(x_i,u_i)$. Then, $y_i \in \mathrm{proj}(\mathcal{I(\mathbb{W})})$. Using the continuity of $f$, we have
	\begin{equation}
	  \begin{split}
	      \mathop {\lim }\limits_{i \to \infty } {y_i} & = \mathop {\lim }\limits_{i \to \infty } f({x_i},{u_i}) = f(\mathop {\lim }\limits_{i \to \infty } {x_i},\mathop {\lim }\limits_{i \to \infty } {u_i}) \\
	      & = f(\bar x,\bar u) \in {\rm{proj(}}{\mathcal{I}}{\rm{(}}\mathbb{W}{\rm{))}},\forall f\in\mathfrak{F}.
	  \end{split}	
	\end{equation}
	Hence, $(\bar x,\bar u) \in{\mathcal{I}}{\rm{(}}\mathbb{W}{\rm{)}}$, which indicates that ${\mathcal{I}}{\rm{(}}\mathbb{W}{\rm{)}}$ is closed. Therefore, ${\mathcal{I}}{\rm{(}}\mathbb{W}{\rm{)}}$ is compact.
	Further more, ${\mathcal{I}^\infty}{\rm{(}}\mathbb{W}{\rm{)}}$ is compact.
	
	To prove that ${\mathcal{I}^\infty}{\rm{(}}\mathbb{W}{\rm{)}}$ is invariant (Definition~\ref{def:invariant}), we need to show that
	\begin{equation}
		\begin{array}{l}
		\forall w \in {{\mathcal I}^\infty }(\mathbb{W}) = \bigcap\limits_{i = 0}^\infty  {{{\mathcal I}^{i + 1}}(\mathbb{W})} ,\\
		f(w) \in {\rm{proj(}}{{\mathcal I}^\infty }(\mathbb{W}){\rm{) = proj(}}\bigcap\limits_{i = 0}^\infty  {{{\mathcal I}^{i + 1}}(\mathbb{W})} {\rm{)}},
		\end{array}
	\end{equation}
	by Definintion~\ref{def:map}, we have for all $i\geq 1$
	\begin{equation}
	\forall w \in {{\mathcal I}^{i + 1}}(\mathbb{W}),f(w,\delta ) \in {\rm{proj(}}{{\mathcal I}^i}(\mathbb{W} ){\rm{)}}\,
	\end{equation}
	which implies
	\begin{equation}
	\forall w \in \bigcap\limits_{i = 0}^\infty  {{{\mathcal I}^{i + 1}}(\mathbb{W})} ,f(w,\delta ) \in \bigcap\limits_{i = 0}^\infty  {{\rm{proj(}}{{\mathcal I}^i}(\mathbb{W}){\rm{)}}}. 
	\end{equation}
	We then aim to show that
	\begin{equation}
		\bigcap\limits_{i = 0}^\infty  {{\rm{proj(}}{{\mathcal I}^i}(\mathbb{W}){\rm{)}}}  = {\rm{proj(}}\bigcap\limits_{i = 0}^\infty  {{{\mathcal I}^i}(\mathbb{W})} {\rm{)}}.
	\end{equation}
	Since ${\rm{proj(}}\bigcap\nolimits_{i = 0}^\infty  {{{\mathcal I}^i}(\mathbb{W})} {\rm{)}} \subset {\rm{proj(}}{{\mathcal I}^i}(\mathbb{W}){\rm{)}}$ for all $i\geq 0 $, we have in general that \[{\rm{proj(}}\bigcap\nolimits_{i = 0}^\infty  {{{\mathcal I}^i}(\mathbb{W})} {\rm{)}} \subset \bigcap\nolimits_{i = 0}^\infty  {{\rm{proj(}}{{\mathcal I}^i}(\mathbb{W}){\rm{)}}} \] and hence it is sufficient to prove the reverse inclusion.
	
	if $x \in \bigcap\nolimits_{i = 0}^\infty  {{\rm{proj}}} ({{\mathcal I}^i}(\mathbb{W}))$, then there exists a sequence \[\mathfrak{U}_i=\{u\in\mathbb{R}^m|\exists x \in {{\rm{proj}}} ({{\mathcal I}^i}(\mathbb{W})), (x,u)\in\mathcal{I}^i(\mathbb{W})\}\] such that $(x,u_j)=w_j\in \mathcal{I}^i(\mathbb{W}), u_j\in\mathfrak{U}_i$. Since $\mathcal{I}(\mathbb{W})$ and the control space $\mathbb{U}$ are compact, by compact definition  in \cite{Kreyszig:1978} (pp.77, 2.5-1), we can find a convergent subsequence $\mathfrak{U}_\infty\subset\mathfrak{U}_i$ and the subsequence $\mathfrak{U}_\infty$ has at least one limit point $\bar u ={\rm{lim}}_{j\rightarrow\infty}u_j$, $\bar u\in\mathfrak{U}_\infty$ such that $(x,\bar u)\in \mathcal{I}^\infty(\mathbb{W})=\bigcap\nolimits_{i = 0}^\infty  ({{\mathcal I}^i}(\mathbb{W}))$, i.e., $x \in {{\rm{proj}}}( \bigcap\nolimits_{i = 0}^\infty  ({{\mathcal I}^i}(\mathbb{W})))$. Therefore,\[\bigcap\nolimits_{i = 0}^\infty  {{\rm{proj(}}{{\mathcal I}^i}(\mathbb{W}){\rm{)}}} \subset {\rm{proj(}}\bigcap\nolimits_{i = 0}^\infty  {{{\mathcal I}^i}(\mathbb{W})} {\rm{)}} ,\] which is what was to be shown.
	\hfill $\blacksquare$	  
\end{pf}

%\begin{myrem}
%	As a topological property, the proof shows that $\mathcal{I}^i(\mathbb{W})$ is compact for all $i\ge 1$ has appeared in~\cite{Rakovic:2006} (p.551, Corollary 2).The relevant proof materials are in Appendix 1~\cite{Rakovic:2006} which is a little cumbersome for the whole proof of this paper. For completeness and clearness, another proof is given in this paper to our perspective.
%\end{myrem}
\begin{myrem}
	To our knowledge, the results similar to Lemma~\ref{lem:I_infty} has been discussed for a half-century, started from the conclusion of the Problem 8 in~\cite{Dugundji:1968}. For infinite-time reachability analysis, a brief proof in state-space is given in Proposition 4 of~\cite{Bertsekas:1972}. The more detailed proof is shown in the Theorem 5.2 of \cite{Franco:2008}(pp. 154). For switched systems, the similar proof is presented in the Theorem 1 of~\cite{Liy:2018}. The result in~\cite{Liy:2018} only requires closedness because the $\mathbb{U}$ of switched systems is finite. Our proof is proposed in a new point of view deals with the subject in state-control space for the purpose of controller design.
\end{myrem}

According to the Lemma~\ref{lem:I_infty}, for a given compact set $\mathbb{W} \subset \mathbb{R}^{n+m}$, $\mathcal{I}^{\infty}(\mathbb{W})$ exists and is the invariant set. However, the negative-definite set $\mathbb{W}_{\mathrm{N}}(L)$ defined in Definition~\ref{def:negative definite} is unbounded and open, namely $\mathbb{W}_{\mathrm{N}}(L)$ is not compact. In order to guarantee the boundedness of $\mathbb{W}_{\mathrm{N}}(L)$, we introduce the following assumption.

\begin{assum} \label{assum:cons_set}
	The state and control input satisfy a set of mixed constrains
	\begin{equation*}
	(x;u) \in \mathbb{W}_\mathrm{cons} \subset \mathbb{R}^{n+m},
	\end{equation*}
	where $\mathbb{W}_\mathrm{cons}$ is a compact set.
\end{assum}

It is without loss of generality to introduce Assumption~\ref{assum:cons_set}, because these constrains typically arise due to physical limitations or safety considerations in practice. Under Assumption~\ref{assum:cons_set}, it follows that $\mathbb{W}_{\mathrm{N}}(L) \subset \mathbb{W}_\mathrm{cons}$, therefore $\mathbb{W}_{\mathrm{N}}(L)$ is bounded. The openness of $\mathbb{W}_{\mathrm{N}}(L)$ is due to that its boundary $\{(x;u)|L(f(x,u)) - L(x) = 0\}$ is not its subset, therefore we modify \eqref{def:w_nl} as 
\begin{eqnarray}
\mathbb{W}_{\rm N}(L)= & \{(x,u) \in \mathbb{R}^{n+m}|\forall x_+ \in\mathbb{X}_+(x,u),\nonumber\\
& L(x_+)-L(x)\le-\alpha \}, \label{def:w_nl<alpha}
\end{eqnarray}

where $\alpha \in \mathbb{R}_+$ is a very small positive constant. Then, it is obvious that, under Assumption~\ref{assum:cons_set}, $\mathbb{W}_{\mathrm{N}}(L)$ defined in \eqref{def:w_nl<alpha} is bounded and closed. Hence, $\mathbb{W}_{\mathrm{N}}(L)$ defined in \eqref{def:w_nl<alpha} is compact. With the compact set $\mathbb{W}_{\mathrm{N}}(L)$, we give the following theorem.

\begin{mytheo}\label{Theo:WNI_Map}
	The $\mathbb{W}_{\mathrm{N\&I}}(L)$ of the state-control space defined in Theorem~\ref{Theo:WNI} satisfies
	\begin{equation} 
	\mathbb{W}_{\mathrm{N\&I}}(L) = \mathcal{I}^\infty(\mathbb{W}_{\mathrm{N}}(L)). \label{eq:W_N&I}
	\end{equation} 
\end{mytheo}

\begin{pf}
	Following Lemma~\ref{lem:I_infty}, we have that the set limit $\mathcal{I}^\infty(\mathbb{W}_{\mathrm{N}}(L)) = \lim_{i \to \infty} \mathcal{I}^{i}(\mathbb{W}_{\mathrm{N}}(L))$ exists and is invariant for plant~\eqref{eq:plant_set}, which means \eqref{eq:ri} is satisfied.
	
	From the definition of the mapping $\mathcal{I}$, we know that
	\begin{equation*}
	\mathcal{I}^\infty(\mathbb{W}_{\mathrm{N}}(L))\subseteq \cdots \subseteq \mathcal{I}^2(\mathbb{W}_{\mathrm{N}}(L)) \subseteq \mathcal{I}(\mathbb{W}_{\mathrm{N}}(L)) \subseteq \mathbb{W}_{\mathrm{N}}(L).
	\end{equation*}
	The above relations mean that $\mathbb{W}_{\mathrm{N\&I}}(L)$ is a subset of $\mathbb{W}_{\mathrm{N}}(L)$. Hence, $\mathcal{I}^{\infty}(\mathbb{W})$ is negative-definite for plant set \eqref{eq:plant_set}, which means \eqref{eq:rnd} is satified. Therefore, we have $\mathbb{W}_{\mathrm{N\&I}}(L) = \mathcal{I}^\infty(\mathbb{W}_{\mathrm{N}}(L))$.
	\hfill $\blacksquare$
	
\end{pf}

Due to the difficulty to obtain the RNIS-SC by its definition, Lemma~\ref{lem:I_infty} and Theorem~\ref{Theo:WNI_Map} are derived to obtain $\mathbb{W}_{\rm N\&I}(L)$ by mapping $\mathcal{I}$. According to the Theorem~\ref{Theo:WNI}, if RNIS-SC is known, the robust controller $\mu$ and the estimate of the closed-loop RDOA ${\rm proj}(\mathbb{W}_{\rm N\&I}(L))$ can be obtained. Then, an interval analysis approach is used to estimate the RNIS-SC, which is introduced in the next subsection. 

\section{Robust Negative-definite and Invariant set Estimation Via Interval Analysis}

In this section, a numerical method estimating the RNIS-SC is proposed. Firstly, an interval analysis approach to estimate a specific set (namely, the SEVIA algorithm), is briefly introduced. Then, an algorithm, estimating the RNIS-SC $\mathbb{W}_{\mathrm{N\&I}}(L)$ for the given Lyapunov function $L$ based on the SEVIA algorithm, is proposed. Thirdly, the convergence of the proposed algorithm is proved under the convergent inclusion function.

\subsection{Set estimation via interval analysis}

Interval analysis is a kind of guaranteed numerical method for approximating sets. Guaranteed means here that approximations of sets of interest are obtained, which can be made as precise as desired \cite{Jaulin:2001}. The basic idea is very simple: use adjacent but disjoint hyperrectangles (rectangles in two-dimensional space, cuboids in three-dimensional space, and hyperrectangles in high-dimensional space) to cover the domain of interest, so as to obtain an approximation of the collection.

The fundamental of interval analysis is the concept of interval vectors and inclusion functions. We briefly introduce these concepts (more details see \cite{Jaulin:2001}).An interval is denoted by
\begin{equation}
{[z_{i}]} = [\underline{z_{i}},\bar z_{i}].\nonumber
\end{equation}
It is a set of real numbers within bounds such that $\underline{z_{i}}\le z_{i}\le \bar {z_{i}}$, where $\underline{z_{i}}$ represents the lower bound ${\rm inf}([z_{i}])$ and the $\bar{z_{i}}$ the upper bound ${\rm sup}([z_{i}])$ of the interval. An interval vector $[z]$ is a subset of $\mathbb{R}^{n_1}$, which is defined as $[z] = [z_{(1)}] \times [z_{(2)}] \times \cdots \times [z_{(n_1)}]$, where the $j$-th interval $[z_{(j)}] = [\underline{z}_{(j)},\bar{z}_{(j)}], j = 1,2,\cdots,n_1$, is a connected subset of $\mathbb{R}$, $\underline{z}_{(j)}$ and $\bar{z}_{(j)}$ are the lower and the upper bound of the interval $[z_{(j)}]$. $[z] \in \mathbb{IR}^{n_1}$ is also called a box, where $\mathbb{IR}^{n_1}$ denotes the set of all $n_1$-dimensional boxes. Considering function $p: \mathbb{R}^{n_1} \to \mathbb{R}^{n_2}$, the interval function $[p]: \mathbb{IR}^{n_1} \to \mathbb{IR}^{n_2}$ is an inclusion function for $p$ if $\forall [z] \in \mathbb{IR}^{n_1}, p([z]) \subset [p]([z])$. An inclusion function $[p]$ is convergent if $\forall [z] \in \mathbb{IR}^{n_1}, \lim_{d([z]) \to 0}d([p]([z])) = 0$, where $d([z]) = \max_{1 \leq j \leq n_1} (\bar{z}_{(j)} - \underline{z}_{(j)})$ denotes the width of box $[z]$. The convergent inclusion function can guarantee the convergence of the set estimation algorithm. Note that for a given function, its convergent inclusion function is not unique, \textit{e.g.}, natural form, centered form and Taylor form.

Suppose the estimated set $\mathbb{W}\subset \mathbb{W}_\mathrm{cons}$ is defined as follows
\begin{eqnarray}
\mathbb{W}=\left\{w\in \mathbb{W}_\mathrm{cons}|w {\quad \rm satisfies\quad CONDITIONS} \right\}.
\end{eqnarray}
The idea of estimating set $\mathbb{W}\subset \mathbb{W}_\mathrm{cons}$ by the SEVIA algorithm is that: first, the CONDITIONS are transformed into the expression of interval inclusion function, so that the given hyper-rectangle $[w] \subset \mathbb{R}^{n+m}$ can be verified whether it entirely in or out $\mathbb{W}$; then, starting from an initial set $\hat{\mathbb{W}}_{\rm init}$ of hyper-rectangles, the SEVIA algorithm performs a recursive exploration, in which four cases may be encountered for the given $[w]$: 
\begin{itemize}
	\item Inner test: if $[w]$ is entirely in $\mathbb{W}$, then $[w]$ is stored in hyperrectangles set $\hat{\mathbb{W}}_{\mathrm{in}}$, as shown in Line 7-8 in Algorithm~\ref{alg:Sevia}.
	\item Outer test: if $[w]$ has an empty intersection with $\mathbb{W}$, then $[z]$ is stored in set $\hat{\mathbb{W}}_{\mathrm{out}}$, as shown in Line 9-10 in Algorithm~\ref{alg:Sevia}. The set collects boxes outside $\mathbb{W}$.
	\item If $[w]$ do not satisfy inner test nor outer test, it means $[w]$ contains the boundary of $\mathbb{W}$; $[w]$ is said to be indeterminate. If the width of $[w]$ is lower than a prespecified parameter $\epsilon > 0$, then it is deemed small enough to be stored in set $\hat{\mathbb{W}}_{\mathrm{bou}}$ collecting boxes containing the boundary of $\mathbb{W}$, as shown in Line 11-12 in Algorithm~\ref{alg:Sevia}.
	\item If $[w]$ is indeterminate and its width is greater than $\epsilon$, then $[w]$ should be bisected and the two newly generated boxes are stored in set $\hat{\mathbb{W}}_\mathrm{do}$ collecting boxes needing further exploration, as shown in Line 13-15 in Algorithm~\ref{alg:Sevia}.
\end{itemize}
The exploration should be recursively implemented until set $\hat{\mathbb{W}}_\mathrm{do}$ is empty.

\begin{algorithm} 
	\caption{Set Estimation Via Interval Analysis} \label{alg:Sevia}
	\begin{algorithmic}[1]
		\Procedure{SEVIA}{Inner test, Outer test, $\hat{\mathbb{W}}_{\mathrm{init}}, \epsilon$}
		\State $\hat{\mathbb{W}}_{\mathrm{in}} := \emptyset, \hat{\mathbb{W}}_{\mathrm{out}} := \emptyset, \hat{\mathbb{W}}_\mathrm{bou} := \emptyset$
		\State $\hat{\mathbb{W}}_{\mathrm{do}} := \hat{\mathbb{W}}_{\mathrm{init}}$
		\While {$\hat{\mathbb{W}}_{\mathrm{do}} \neq \emptyset$}
		\State Get a box $[w]$ from $\hat{\mathbb{W}}_{\mathrm{do}}$
		\State Remove $[w]$ from $\hat{\mathbb{W}}_{\mathrm{do}}$
		\If {$[w]$ satisfies Inner test}
		\State Add $[w]$ to set $\hat{\mathbb{W}}_{\mathrm{in}}$
		\ElsIf {$[w]$} satisfies Outer test
		\State Add $[w]$ to set $\hat{\mathbb{W}}_{\mathrm{out}}$
		\ElsIf {$d([w]) < \epsilon$}
		\State Add $[w]$ to set $\hat{\mathbb{W}}_{\mathrm{bou}}$
		\Else
		\State Bisect box $[w]$
		\State Add the two new boxes to set $\hat{\mathbb{W}}_{\mathrm{do}}$
		\EndIf
		\EndWhile
		\State \textbf{return} $\hat{\mathbb{W}}_{\mathrm{in}}$ as the estimate of $\mathbb{W}$%, Z_{\mathrm{bou}}, Z_{\mathrm{out}}$
		\EndProcedure
	\end{algorithmic}
\end{algorithm}

\subsection{Robust negative-definite and invariant set estimation via interval analysis }

In this subsection, the RNISEVIA (Robust negative-definite and invariant set estimation via interval analysis) algorithm is proposed to estimate the robust negative-definite and invariant set $\mathbb{W}_{\mathrm{N\&I}}(L)$ for the given Lyapunov function $L$ based on the SEVIA algorithm. 

The inclusion function of the function defined from Euclidean space to Euclidean space is well defined and researched \cite{Jaulin:2001}. However, the definition of the inclusion function of $\mathbb{X}_+:\mathbb{R}^{n+m}\to\mathbb{IR}^n$ defined in \eqref{eq:nextset_x} does not exist in literatures. According to the definition of $\mathbb{X}_+(w)$ in~\eqref{eq:nextset_x}, we define the interval function $[\mathbb{X}_+]:\mathbb{IR}^{n+m}\to\mathbb{IR}^n$ as
\begin{eqnarray}
	[\mathbb{X}_+]([w])=[\underline{f}]([w])\sqcup [\bar{f}]([w]), \label{defn:[X]_+}
\end{eqnarray}
where $\underline{f}:\mathbb{R}^{n+m} \to \mathbb{R}^n$ is defined as $\underline{f}(w)=\hat{f}(w)-\delta(w)$ and $\bar{f}:\mathbb{R}^{n+m} \to \mathbb{R}^n$ is defined as $\bar{f}(w)=\hat{f}(w)+\delta(w)$. Actually, the interval function $[\mathbb{X}_+]$ is the inclusion function of $\mathbb{X}_+$ as shown in the following lemma.

\begin{mylem}\label{Lem:[X]}
	Interval function $[\mathbb{X}_+]:\mathbb{IR}^{n+m}\rightarrow\mathbb{IR}^n$ defined in \eqref{defn:[X]_+} is the inclusion function of $\mathbb{X}_+:\mathbb{R}^{n+m}\rightarrow\mathbb{R}^n$ defined in \eqref{eq:nextset_x}, namely, $\forall [w]\in\mathbb{IR}^{n+m},\mathbb{X}_+([w])\subset[\mathbb{X}_+]([w])$.
\end{mylem}

\begin{pf}
	According to the~\eqref{eq:nextset_x}, we know that $\mathbb{X}_+(w)=[\underline{f}(w),\bar{f}(w)]$, hence we have 	
	\begin{equation}
	\mathbb{X}_+([w])=[{\rm inf}(\underline{f}([w])),{\rm sup}(\bar f([w]))]. \nonumber \label{eq:def_X_[w]}
	\end{equation}
	According to the properties of inclusion function, we know that
	\begin{eqnarray}
	\forall [w]\in\mathbb{IR}^{n+m},\underline{f}([w])\subset[\underline{f}]([w]),\nonumber\\\forall [w]\in\mathbb{IR}^{n+m},\bar{f}([w])\subset\bar{[f]}([w]).\nonumber
	\end{eqnarray}
	Therefore, 
	\begin{eqnarray}
	&\forall [w]\in\mathbb{IR}^{n+m}, {\rm inf}(\underline{f}([w]))\in[\underline{f}]([w])\nonumber,\\&\forall [w]\in\mathbb{IR}^{n+m},{\rm sup}(\bar{f}([w]))\in\bar{[f]}([w]).\nonumber
	\end{eqnarray}
	which means
	\begin{eqnarray}
	[{\rm inf}(\underline{f}([w])),{\rm sup}(\bar f([w]))]\subset[\underline{f}]([w])\sqcup [\bar{f}]([w]). \nonumber
	\end{eqnarray}
	Hence, we have $\forall [w]\in\mathbb{IR}^{n+m},\mathbb{X}_+([w])\subset[\mathbb{X}_+]([w])$. \hfill $\blacksquare$	
\end{pf}

Once the inclusion function $[\mathbb{X}_+]([w])$ of $\mathbb{X}_+(w)$ is given, the Inner test and Outer test in SEVIA can be obtained. For a given Lyapunov function $L$, the RNS-SC $\mathbb{W}_{\rm N}(L) \subset \mathbb{R}^n$ can be approximated using SEVIA algorithm, as shown in Line 2-3 of Algorithm \ref{alg:est_W_{NI}_(L)}. By~\eqref{def:w_nl<alpha}, if $w=(x,u)\in\mathbb{W}_{\rm N}(L)$, $w$ satisfies
\begin{eqnarray}
L(\mathbb{X}_+(w))-L(x)\subset(-\infty,-\alpha],
\end{eqnarray}
where $L(\mathbb{X}_+(w))$ denotes the range of the function $L$ whose set is defined by $\mathbb{X}_+(w)$. Hence, if $[w]=[x]\times[u]\subset\mathbb{W}_{\rm N}(L)$, then $[w]$ must satisfy
\begin{eqnarray}
L(\mathbb{X}_+([w]))-L([x])\subset(-\infty,-\alpha]. \label{eq:L(X_+)-L(x)}
\end{eqnarray}

Depending on the properties of the interval operation and the inclusion function, we can get
\begin{eqnarray}
L(\mathbb{X}_+([w]))-L([x])\subset[L]([\mathbb{X}_+]([w]))-[L]([x]). \label{eq:[L]_subset}
\end{eqnarray}
By~\eqref{eq:[L]_subset}, if
\begin{eqnarray}
[L]([\mathbb{X}_+]([w]))-[L]([x])\subset(-\infty,-\alpha], \label{eq:[L]<0}
\end{eqnarray}
then~\eqref{eq:L(X_+)-L(x)} must be satisfied so that $[w]\subset\mathbb{W}_{\rm N}(L)$. Similarly, if
\begin{eqnarray}
[L]([\mathbb{X}_+]([w]))-[L]([x])\cap(-\infty,-\alpha]=\emptyset, \label{eq:[L]cap0}
\end{eqnarray}
then $L(\mathbb{X}_+([w]))-L([x])$ must be outside $(-\infty,-\alpha]$ so that $[w]$ is outside $\mathbb{W}_{\rm N}(L)$. With the internal validation expressions~\eqref{eq:[L]<0} and external validation expressions~\eqref{eq:[L]cap0}, the estimated $\hat{\mathbb{W}}_{\rm N}(L)$ of RNS-SC $\mathbb{W}_{\rm N}(L)$ can be derived from the following formula:
\begin{eqnarray}
\hat{\mathbb{W}}_{\rm N}(L):={\rm SEVIA}(\eqref{eq:[L]<0},\eqref{eq:[L]cap0},\{[w_{\rm init}]\},\epsilon), \label{eq:sevia_W_n}
\end{eqnarray}
where $[w_{\rm init}]\subset\mathbb{R}^{n+m}$ is a hyperrectangle large enough in the state-control space.

With the initial set $\hat{\mathbb{W}}_\mathrm{N}(L)$, recursively using SEVIA algorithm to approximate mapping $\mathcal{I}$, an inner approximation $\hat{\mathbb{W}}_{\mathrm{N\&I}}(L)$ of the negative-definite and invariant set $\mathbb{W}_{\mathrm{N\&I}}(L)$ defined in \eqref{eq:ri} can be obtained, as shown in Line 4-10 of Algorithm \ref{alg:est_W_{NI}_(L)}. For any hyperrectangular set $\hat{\mathbb{W}}_\mathrm{N}(L)\subset \mathbb{R}^{n+m}$, we know from \eqref{eq:defn:I} that if $[w]\subset\mathcal{I(\hat{\mathbb{W}})}_1$, $[w]\subset\hat{\mathbb{W}}_1$ needs to satisfy 
\begin{eqnarray}
\mathbb{X}_+([w])\subset{\rm proj}(\hat{\mathbb{W}}_1). \label{eq:x_+subproj}
\end{eqnarray}
Due to the difficulty of solving the function range, the verification of the above formula is hard to achieve. Therefore, the interval inclusion function is used instead of the solution of the function range. Based on $\mathbb{X}_+([w])\subset[\mathbb{X}_+]([w])$, we know that if
\begin{eqnarray}
[\mathbb{X}_+]([w])\subset{\rm proj}(\hat{\mathbb{W}}_1), \label{eq:[x]_+subproj}
\end{eqnarray}
then \eqref{eq:x_+subproj} must be satisfied, therefore $[w]\subset\mathcal{I}(\hat{\mathbb{W}}_1)$. Similarly, if
\begin{eqnarray}
[\mathbb{X}_+]([w])\cap{\rm proj}(\hat{\mathbb{W}}_1)=\emptyset, \label{eq:[x]_capproj}
\end{eqnarray}
then $[\mathbb{X}_+]([w])$ must be outside ${\rm proj}(\hat{\mathbb{W}}_1)$, therefore $[w]$ is outside $\mathcal{I}(\hat{\mathbb{W}}_1)$. With internal validation expressions~\eqref{eq:[x]_+subproj} and external validation expressions~\eqref{eq:[x]_capproj}, the estimated $\hat{\mathbb{W}}_2$ of $\mathcal{I}(\hat{\mathbb{W}}_1)$ can be derived from the following formula:
\begin{eqnarray}
\hat{\mathbb{W}}_2:={\rm SEVIA}(\eqref{eq:[x]_+subproj},\eqref{eq:[x]_capproj},\hat{\mathbb{W}}_1,\epsilon)
\end{eqnarray}

\begin{algorithm} 
	\caption{Robust Negative-definite and Invariant Set Estimation Via Interval Analysis \label{alg:est_W_{NI}_(L)}}
	\begin{algorithmic}[1]
		\Procedure{RNISEVIA}{$\mathfrak{F}, L, [w_{\mathrm{init}}], \epsilon$}
		\State Using $\mathfrak{F}$ and $L$ to get the internal validation expressions~\eqref{eq:[L]<0} and external validation expressions~\eqref{eq:[L]cap0}.
		\State $\hat{\mathbb{W}}_{\rm N}(L):={\rm SEVIA}(\eqref{eq:[L]<0},\eqref{eq:[L]cap0},\{[w_{\rm init}]\},\epsilon)$
		\State Using $\mathfrak{F}$ and $L$ to get the internal validation expressions~\eqref{eq:[x]_+subproj} and external validation expressions~\eqref{eq:[x]_capproj}.
		\State $\hat{\mathbb{W}}_1 := \emptyset$
		\State $\hat{\mathbb{W}}_2 := \hat{\mathbb{W}}_{\rm N}(L)$
		\While{$\hat{\mathbb{W}}_1 \neq \hat{\mathbb{W}}_2$}
		\State $\hat{\mathbb{W}}_1 := \hat{\mathbb{W}}_2$
		\State $\hat{\mathbb{W}}_2 := {\rm SEVIA}(\eqref{eq:[x]_+subproj},\eqref{eq:[x]_capproj},\hat{\mathbb{W}}_1,\epsilon)$	
		\EndWhile
		\State $\hat{\mathbb{W}}_{\mathrm{N\&I}}(L) := \hat{\mathbb{W}}_1$
		\State \textbf{return} $\hat{\mathbb{W}}_{\mathrm{N\&I}}(L)$ as the internal approximation of $\mathbb{W}_{\rm N\&I}(L)$
		\EndProcedure
	\end{algorithmic}
\end{algorithm}

\subsection{Convergence of RNISEVIA algorithm}

RNIDEVIA is a finite algorithm, which terminates in less than $(d([w_{init}])/\epsilon)^2\times(d([w_{init}])/\epsilon+1)$ times. Let us prove the most important properties of RNIDEVIA: the convergence of the output result $\hat{\mathbb{W}}_{\mathrm{N\&I}}(L)$.

According to the~\eqref{defn:[X]_+}, $[\mathbb{X}_+]([w])$ is the interval union of two inclusion functions in Euclidean space. Three forms (natural, center and Taylor) of inclusion function are convergent~\cite{Jaulin:2001} (pp. 27-38). For the interval union $[g]\sqcup[h]$ of two inclusion fuctions $[g]$ and $[h]$, this propertiy can be extended as follows:

\begin{mydefn}  \label{Defn:convergent}
	An inclusion function $[\mathbb{X}_+]([w])=[\underline{f}]([w])\sqcup [\bar{f}]([w]):\mathbb{IR}^{n+m}\rightarrow\mathbb{IR}^n$ for $[\underline{f}]([w]):\mathbb{IR}^{n+m}\rightarrow\mathbb{IR}^n$ and $[\bar{f}]([w]):\mathbb{IR}^{n+m}\rightarrow\mathbb{IR}^n$ is convergent if, for any sequence of boxes $[w](j)$,
	\begin{equation}
	\mathop {\lim }\limits_{j \to \infty } d([w](j)) = 0 \Rightarrow \mathop {\lim }\limits_{j \to \infty } d([\mathbb{X}_+](j)) = D, \label{eq:convergent}
	\end{equation}
	where $d([w])$ is the width of $[w]$, $D={\rm sup}(\bar f([w]))-{\rm inf}(\underline{f}([w]))$.
\end{mydefn}  

\begin{mylem} \label{Lem:[X]convergence}
	If the inclusion functions $[\underline{f}]([w]), [\bar f]([w])$ are in natural, center or Taylor form, the inclusion function $[\mathbb{X}_+]([w])$ satisfies~\eqref{defn:[X]_+} is convergent.
\end{mylem}

\begin{pf}
	
	We have $[\mathbb{X}_+]([w])=[\underline{f}]([w])\sqcup [\bar{f}]([w])$, hence, for any sequence of boxes $[w](j)$
	\begin{eqnarray}
	\mathop {\lim }\limits_{j \to \infty } d([w](j)) = 0 \Rightarrow \mathop {\lim }\limits_{j \to \infty } d([\mathbb{X}_+]([w](j)))\nonumber
	\\ = sup([\bar f]([w](j)))-inf([\underline{f}]([w](j))),  \label{Pf:[X]}
	\end{eqnarray}
	where $d([x])$ is the width of $[x]$. Since the inclusion functions $[\underline{f}]([w]), [\bar f]([w])$ are in natural, center or Taylor form, then $[\underline{f}], [\bar f]$ are convergent which indicate that, 
	\begin{eqnarray}
	&\mathop {\lim }\limits_{j \to \infty } d([w](j)) = 0 \Rightarrow \mathop {\lim }\limits_{j \to \infty } d([\underline{f}]([w(j)])) = 0,\nonumber
	\\&\mathop {\lim }\limits_{j \to \infty } d([w](j)) = 0 \Rightarrow \mathop {\lim }\limits_{j \to \infty } d([\bar f]([w](j))) = 0.\nonumber
	\end{eqnarray}
	Therefore,
	\begin{eqnarray}
	&\mathop {\lim }\limits_{j \to \infty } d([w](j)) = 0 \Rightarrow \mathop {\lim }\limits_{j \to \infty } [\underline{f}]([w](j)) = \underline{f}(w) ,\nonumber
	\\&\mathop {\lim }\limits_{j \to \infty } d([w](j)) = 0 \Rightarrow \mathop {\lim }\limits_{j \to \infty } [\bar{f}]([w](j)) = \bar{f}(w).\nonumber
	\end{eqnarray}
	Let $D={\rm sup}(\bar f([w]))-{\rm inf}(\underline{f}([w]))$, then~\eqref{Pf:[X]} can be rewritten as
	\begin{eqnarray}
	\mathop {\lim }\limits_{j \to \infty } d([w](j)) = 0 \Rightarrow \mathop {\lim }\limits_{j \to \infty } d([\mathbb{X}_+]([w](j))) = D,\nonumber
	\end{eqnarray}
	which is what was to be shown.
	\hfill $\blacksquare$	
\end{pf}

\begin{mydefn}(Definition 13 in \cite{Jaulin:1993})
	The compact $\mathbb{A}$ is full if ${\rm clo(int(\mathbb{A}))}=\mathbb{A}$, where ${\rm clo(\mathbb{A})}$ is closure of $\mathbb{A}$ and ${\rm int(\mathbb{A})}$ is interior of the set $\mathbb{A}$.
\end{mydefn}

\begin{mytheo} \label{Theo:inclusionfunction}
	If $\mathbb{W}_{\rm N\&I}(L)$ is full, then the set $\hat{\mathbb{W}}_{\mathrm{N\&I}}(L)$ evaluated by RNIDEVIA ($\mathfrak{F}, L, [w_{\mathrm{init}}], \epsilon$) satisfies the following property:
	\begin{equation}
	\mathop {\lim }\limits_{\epsilon  \to 0}  \hat{\mathbb{W}}_{\mathrm{N\&I}}(L)=\mathbb{W}_{\rm N\&I}(L). \nonumber
	\end{equation}
\end{mytheo}

\begin{pf}
	First we show that the $\hat{\mathbb{W}}_{\rm N}(L)$ obtained at Line 3 converges to $\mathbb{W}_{\rm N}(L)$ when $\epsilon\rightarrow 0$. There are two funcitons in SEVIA of~\eqref{eq:sevia_W_n}: $[L]$ and $[\mathbb{X}_+]$. The convergence of $[\mathbb{X}_+]([w]):\mathbb{IR}^{n+m}\rightarrow\mathbb{IR}^n$ is proved in Lemma~\ref{Lem:[X]convergence} such that $
	\mathbb{X}_+([w])\subset[\mathbb{X}_+]([w])$ and $ \lim_{d([w]) \to 0}[\mathbb{X}_+]([w])=\mathbb{X}_+([w])$. For a Lyapunov function $L$ in Euclidean space, its interval inclusion function $[L]:\mathbb{IR}^{n}\rightarrow\mathbb{IR}$ satisfies that $[L](\mathbb{X}_+([w]))\subset [L]([\mathbb{X}_+]([w]))$ and $ \lim_{[\mathbb{X}_+]([w])\to\mathbb{X}_+([w])}[L]([\mathbb{X}_+]([w]))= [L](\mathbb{X}_+([w]))$. Therefore, we can get that
	$\lim_{d([w])  \to 0}$ $  [L]([\mathbb{X}_+]([w]))= L(\mathbb{X}_+([w])),\lim_{d([w])  \to 0}$ $ [L]([x])=L([x])$.
	
	Once the SEVIA algorithm at Line 3 is finished, if $[w]\in\hat{\mathbb{W}}_{\rm bou}$, then $d([w])<\epsilon$. When $\epsilon\rightarrow 0$, we have $d([w])\rightarrow 0$. The inner test~\eqref{eq:[L]<0} and outer test~\eqref{eq:[L]cap0} can be rewritten as
	\begin{eqnarray}
	L(\mathbb{X}_+([w]))-L([x])\subset(-\infty,-\alpha], \label{eq:inner} \\
	L(\mathbb{X}_+([w]))-L([x])\cap(-\infty,-\alpha]=\emptyset. \label{eq:outer}
	\end{eqnarray}
	Since the $[w]\in\hat{\mathbb{W}}_{\rm bou}$  is not satisfied with \eqref{eq:inner} or \eqref{eq:outer}, when $\epsilon\rightarrow 0$, then $L(\mathbb{X}_+([w]))-L([x])=0$.
	Consequently,
	\begin{eqnarray}
	\forall[w]\in\hat{\mathbb{W}}_{\rm bou},& \lim_{\epsilon\rightarrow 0}h_\infty([w],\partial\mathbb{W}_{\rm N}(L))= 0 \nonumber\\\Rightarrow 
	 &\lim_{\epsilon\rightarrow 0}h_\infty(\hat{\mathbb{W}}_{\mathrm{bou}},\partial\mathbb{W}_{\rm N}(L))= 0 ,
	\end{eqnarray}
	where $\partial\mathbb{W}_{\rm N}(L)$ denotes the boundary of the compact set $\mathbb{W}_{\rm N}(L))$, $h_\infty(\hat{\mathbb{W}}_{\mathrm{bou}},\partial\mathbb{W}_{\rm N}(L))$ represents the Hausdorff distance (\textit{e.g.}, Definition 10 in~\cite{Jaulin:1993}) between $\hat{\mathbb{W}}_{\mathrm{bou}}$ and $\partial\mathbb{W}_{\rm N}(L))$. The Hausdorff distance based on the infinity norm is metric for the set of compact subsets of $\mathbb{R}^{n+m}$. It is trival that $\partial\mathbb{W}_{\rm N}(L)\subset\hat{\mathbb{W}}_{\mathrm{bou}}$. If $\epsilon = 0$ then $h_\infty(\hat{\mathbb{W}}_{\mathrm{bou}},\partial\mathbb{W}_{\rm N}(L))=0$ thus $h_\infty(\hat{\mathbb{W}}_{\mathrm{bou}}=\partial\mathbb{W}_{\rm N}(L))$, which gives $\lim_{\epsilon  \to 0}\hat{\mathbb{W}}_{\mathrm{bou}}=\partial\mathbb{W}_{\rm N}(L))$. 
	
	Now that the boundary $\partial\mathbb{W}_{\rm N}(L)$ has been determined, the rest of the boxes are either outside~\eqref{eq:inner} or inside~\eqref{eq:outer}. Note that the $\mathbb{W}_{\rm N\&I}(L)$ is full and $\mathbb{W}_{\rm N\&I}(L)\subset\mathbb{W}_{\rm N}$. Hence the $\mathbb{W}_{\rm N}$ is also full which indicates that there exists finite boxes $\hat{\mathbb{W}}_{\mathrm{in}}$ satisfies \eqref{eq:inner}. Since $\lim_{\epsilon  \to 0}\hat{\mathbb{W}}_{\mathrm{bou}}=\partial\mathbb{W}_{\rm N}(L)$, the rest boxes are either inner or outer. We can finally get $\lim_{\epsilon  \to 0}\hat{\mathbb{W}}_{\rm N}(L) := \hat{\mathbb{W}}_{\mathrm{in}}=\mathbb{W}_{\rm N}(L)$. Hence, with convergent inclusion function $[\mathbb{X}_+]([w])$, SEVIA can ganrantee the convergence of the output when $\epsilon\rightarrow0$.
	
	Similarly, no matter how many times the SEVIA algorithm executed in the loop in Line 7-10 of Algorithm~\ref{alg:est_W_{NI}_(L)}, since the inclusion function $[\mathbb{X}_+]([w])$ in~\eqref{eq:[x]_capproj} is convergent, each SEVIA satisfies
	\begin{equation}
	\lim_{\epsilon\rightarrow 0}h_\infty(\hat{\mathbb{W}}_{\mathrm{bou}},\partial \mathbb{W}) = 0, \nonumber
	\end{equation}
	which indicates that $h_\infty(\hat{\mathbb{W}}_2, \mathbb{W}_2)$ tends to zero when $\epsilon = 0$. Hence, the final output $\hat{\mathbb{W}}_{\mathrm{N\&I}}(L) := \hat{\mathbb{W}}_1= \hat{\mathbb{W}}_2$ satisfies
	\begin{equation}
	\lim_{\epsilon\rightarrow 0}\  h_\infty(\hat{\mathbb{W}}_{\mathrm{N\&I}}(L),\mathbb{W}_{\rm N\&I}(L))= 0.\nonumber
	\end{equation}
	As an inner approximation of full compact $\mathbb{W}_{\rm N\&I}(L)$, we have $\hat{\mathbb{W}}_{\mathrm{N\&I}}(L)\subset\mathbb{W}_{\rm N\&I}(L)$. If $\epsilon = 0$ then $h_\infty (\hat{\mathbb{W}}_{\mathrm{N\&I}}(L),$ $\mathbb{W}_{\rm N\&I}(L)) =0$ thus $\hat{\mathbb{W}}_{\mathrm{N\&I}}(L)=\mathbb{W}_{\rm N\&I}(L)$. Consequently, $\lim_{\epsilon  \to 0}  \hat{\mathbb{W}}_{\mathrm{N\&I}}(L)=\mathbb{W}_{\rm N\&I}(L)$.
	\hfill $\blacksquare$
	
\end{pf}

Depending to the Theorem~\ref{Theo:inclusionfunction}, the inner approximate $\hat{\mathbb{W}}_{\mathrm{N\&I}}(L)$ converges to set $\mathbb{W}_{\rm N\&I}(L)$ when $\epsilon$ tends to zero. This means that RNIDEVIA can approximate RNIS-SC $\mathbb{W}_{\rm N\&I}(L)$ with an arbitrary precision similar to the cyclotomic method. In practice, note that there is of course a limit to the accuracy that can be achieved, because of the complexity of the resulting description \cite{Kieffer:2002}. 

With the RNIS-SC estimate method, the RDOA for a given Lyapunov function can be obtained. Hence, the RDOA enlargement approach is introduced in the next section.

\section{Closed-loop RDOA Enlargement and Controller Design}

\subsection{Closed-loop RDOA Enlargement}

According to the Theorem~\ref{Theo:WNI}, for a given Lyapunov function $L$, RNID-SC $\mathbb{W}_{\mathrm{N\&I}}(L)$ is a robust stablization controller set. The projection ${\rm proj}(\mathbb{W}_{\mathrm{N\&I}}(L))$ of $\mathbb{W}_{\mathrm{N\&I}}(L)$ along the state space is an estimate of RDOA. The estimate $\hat{\mathbb{W}}_{\rm N\&I}(L)$ of $\mathbb{W}_{\mathrm{N\&I}}(L)$ can be obtained by using Algorithm~\eqref{alg:est_W_{NI}_(L)}, then the estimate of closed-loop RDOA is ${\rm proj} (\hat{\mathbb{W}}_{\rm N\&I}(L))$. 

Inspired by \cite{Li:2020}, we observed that, for different Lyapunov functions, the robust negative-definite and invariant set  $\hat{\mathbb{W}}_{\rm N\&I}(L)$ are totally different so as its projection ${\rm proj} (\hat{\mathbb{W}}_{\rm N\&I}(L))$. Therefore, the motivation of this subsection is to enlarge the estimate of the closed-loop RDOA by selecting an appropriate Lyapunov function. In order to achieve this goal, a solvable optimization problem is formulated to select an appropriate Lyapunov function from a parameterized positive-definite function set, which is a subset of all sum-of-square polynomials~\cite{Powershot:1998}. The objective function of this optimization problem is the volume of the estimate of the closed-loop RDOA, which can be obtained based on the method proposed in the last subsection. The analytical expression of the volume of the estimate of the closed-loop RDOA is hard to be derived, but it is easy to evaluate for a given positive-definite function. Hence, if a positive-definite function set rather than a Lyapunov function is given, a significantly large estimate of the closed-loop RDOA may be obtained by selecting an appropriate Lyapunov function from the positive-definite function set~\cite{Topcu:2008}. Based on this idea, the following optimization problem is formulated.

\begin{equation}
\max_{L \in \mathfrak{L}_{n,2d}} \mathfrak{m}(\mathrm{proj}(\hat{\mathbb{W}}_{\mathrm{N\&I}}(L))), \label{eq:optim}
\end{equation}
where $\mathfrak{m}(\mathrm{proj}(\hat{\mathbb{W}}_{\mathrm{N\&I}}(L)))$ denotes the Lebesgue measure of $\mathrm{proj}(\hat{\mathbb{W}}_{\mathrm{N\&I}}(L))$ (in Euclidean space, it is the volume) and $\mathfrak{L}_{n,2d}$ is a subset of all sum-of-square polynomials \cite{Ufuk:2010} in $n$ variables with degree $\leq 2d$, which is defined as
\begin{equation}
\mathfrak{L}_{n,2d} = \Big\{L \in \mathfrak{R}_{n,2d} \Big| L(x) = S^T_d(x)P^TPS_d(x), x \in \mathbb{R}^n\Big\}, \nonumber
\end{equation}
where $\mathfrak{R}_{n,2d}$ denotes the set of all polynomials in $n$ variables with degree $\leq 2d$, $P \in \mathbb{R}^{r \times r}$ has full rank,
\begin{equation}
S_d(x) = (x_{(1)};\cdots;x_{(n)};x_{(1)} x_{(2)};\cdots;x^d_{(n)}) \in \mathbb{R}^r, \nonumber
\end{equation}
and $r = \left(
\begin{smallmatrix}
\scriptscriptstyle n+d \\
\scriptscriptstyle d
\end{smallmatrix} \right) - 1$ represents the dimension of $P$.

Two points should be noted about the optimization problem~\eqref{eq:optim}. First, for the given Lyapunov function $L$, the volume of $\mathrm{proj}(\hat{\mathbb{W}}_{\mathrm{N\&I}}(L))$ is easy to calculate, since $\mathrm{proj}(\hat{\mathbb{W}}_{\mathrm{N\&I}}(L))$ consists of boxes. Second, the Lyapunov function $L$ is selected from a parameterized polynomial function set with the parameters $P \in \mathbb{R}^{r \times r}$.

Since $L(x) = S^T_d(x)P^TPS_d(x)$ and $\hat{\mathbb{W}}_{\mathrm{N\&I}}(L)$ can be obtained by Algorithm~\ref{alg:est_W_{NI}_(L)}, $P$ is our major concern. We defined function $m: \mathbb{R}^{r \times r} \to \mathbb{R}$ as 
\begin{eqnarray}
   m(P) = \mathfrak{m}(\mathrm{proj}(\hat{\mathbb{W}}_{\mathrm{N\&I}}(L))).\label{eq:mp}
\end{eqnarray}
With the function~\eqref{eq:mp}, the optimization problem~\eqref{eq:optim} can be equivalently rewritten as
\begin{equation}
\max_{P \in \mathbb{R}^{r \times r}} m(P). \label{eq:optim_P}
\end{equation}
The analytical expression of $m(P)$ is hard to be derived, but it is easy to evaluate $m(P)$ for a given $P$. Hence, classic optimization methods, e.g., gradient descent method, cannot be used to solve the optimization problem~\eqref{eq:optim_P}. However, meta-heuristic optimization methods can be used to solve the optimization problem~\eqref{eq:optim_P}, whose advantage is that the function to be optimized is only required to be evaluable. Popular meta-heuristic optimizers for real-valued search-spaces include particle swarm optimization~\cite{Kennedy:1995}, differential evolution~\cite{Opara:2019} and evolution strategies\cite{Wang:2011}. There are lots of literatures about meta-heuristic optimizers\cite{Marler:2004}, so we omit the introduction about them in this paper.

\subsection{Controller Design}

Since $L(f(0,0)) - L(0) = 0$, the origin $0 \in \mathbb{R}^{n+m}$ is in the boundary of $\mathbb{W}_{\mathrm{N\&I}}(L)$. Hence, there is no box belonging to the inner approximation $\hat{\mathbb{W}}_{\mathrm{N\&I}}(L)$ nearby the origin $0 \in \mathbb{R}^{n+m}$ and there is a small neighborhood $\mathbb{X}_0$ of the origin $0 \in \mathbb{R}^{n}$ that is not contained by $\mathrm{proj}(\hat{\mathbb{W}}_{\mathrm{N\&I}}(L))$. The size of the neighborhood $\mathbb{X}_0$ depends on parameter $\epsilon$ in Algorithm~\ref{alg:est_W_{NI}_(L)}. Since we suppose that the linearization of plant set \eqref{eq:plant_set} is controllable at the origin, there must exist a linear controller which is able to stabilize all state in $\mathbb{X}_0$ when the size of $\mathbb{X}_0$ is small enough. As a result, the estimate of the closed-loop DOA should be $\mathrm{proj}(\hat{\mathbb{W}}_{\mathrm{N\&I}}(L)) \cup \mathbb{X}_0$ and the corresponding state feedback controller is 
\begin{equation} \label{eq:controller}
\mu(x) = \left\{ 
\begin{array}{ll}
Kx, & \mathrm{if \ } x \in \mathbb{X}_0 \\
\tilde{\mu}(x), & \mathrm{if \ } x \in \mathrm{proj}(\hat{\mathbb{W}}_{\mathrm{N\&I}}(L))
\end{array}
\right.,
\end{equation}
where $K \in \mathbb{R}^{m \times n}$ is obtained by the linear controller design method according to the linearization of plant \eqref{eq:plant_set} and $\tilde{\mu}: \mathbb{R}^n \to \mathbb{R}^m$ satisfies $\forall x \in \mathrm{proj}(\hat{\mathbb{W}}_{\mathrm{N\&I}}(L)), (x;\tilde{\mu}(x)) \in \hat{\mathbb{W}}_{\mathrm{N\&I}}(L)$. 

According to the Theorem~\ref{Theo:WNI}, all controllers belonging to $\hat{\mathbb{W}}_{\mathrm{N\&I}}(L)$ can stabilize the plant. In order to verify the effectiveness of our method, the controller design method is proposed as follows. 

First, select a training set belonging to $\hat{\mathbb{W}}_{\mathrm{N\&I}}(L)$; then, obtain $\mu$ with a function estimation method, such as interpolation, Gaussian processes regression and so on. When the trend of the training data points is smooth and $\mu(0)=0$, it can be guaranteed that $\tilde{\mu}$ obtained from the function estimator satisfies $\forall x \in \mathrm{proj}(\hat{\mathbb{W}}_{\mathrm{N\&I}}(L)), (x;\tilde{\mu}(x)) \in \hat{\mathbb{W}}_{\mathrm{N\&I}}(L)$. 

Note that the problem of finding an optimal controller is not in the scope of this paper. Based on the controller design method, the simulation of a specific plant is formulated in the next section.

\section{Simulation}

Consider the nominal model $\hat{f}$ and the modeling error bound $\delta(x,u)$ of plant set~\eqref{eq:plant_set} is:
\begin{eqnarray}
\hat{f}(x,u) &=& -\sin(2x) - xu - 0.2x - u^2 + u, \nonumber \\
\delta(x,u) &=& 1 - \exp \left(-0.5(x^2 + u^2)\right), \nonumber
\end{eqnarray} 
where $x \in \mathbb{R}$ and $u \in \mathbb{R}$. The interested region is $\mathbb{W}=[-2,2]\times[-2,2]\subset\mathbb{R}^2$ in the state-control space.

The nominal model is the same as~\cite{Li:2014_79, Li:2019} for the study on the effect of modeling error, the modeling error bound is the same as~\cite{Li:2020} for result comparison. To reduce the conservatism of interval arithmetic, the intervals are split, which in turn results in an exponential complexity with respect to the number of state variables. For reference to other systems, the computing times are provided in Appendix~\ref{Appendix} with pseudocode. 

\subsection{Stabilization with closed-loop DOA estimation for a given Lyapunov function}
\begin{figure*}[ht]
	\begin{center}
		\includegraphics[width=0.295\textwidth]{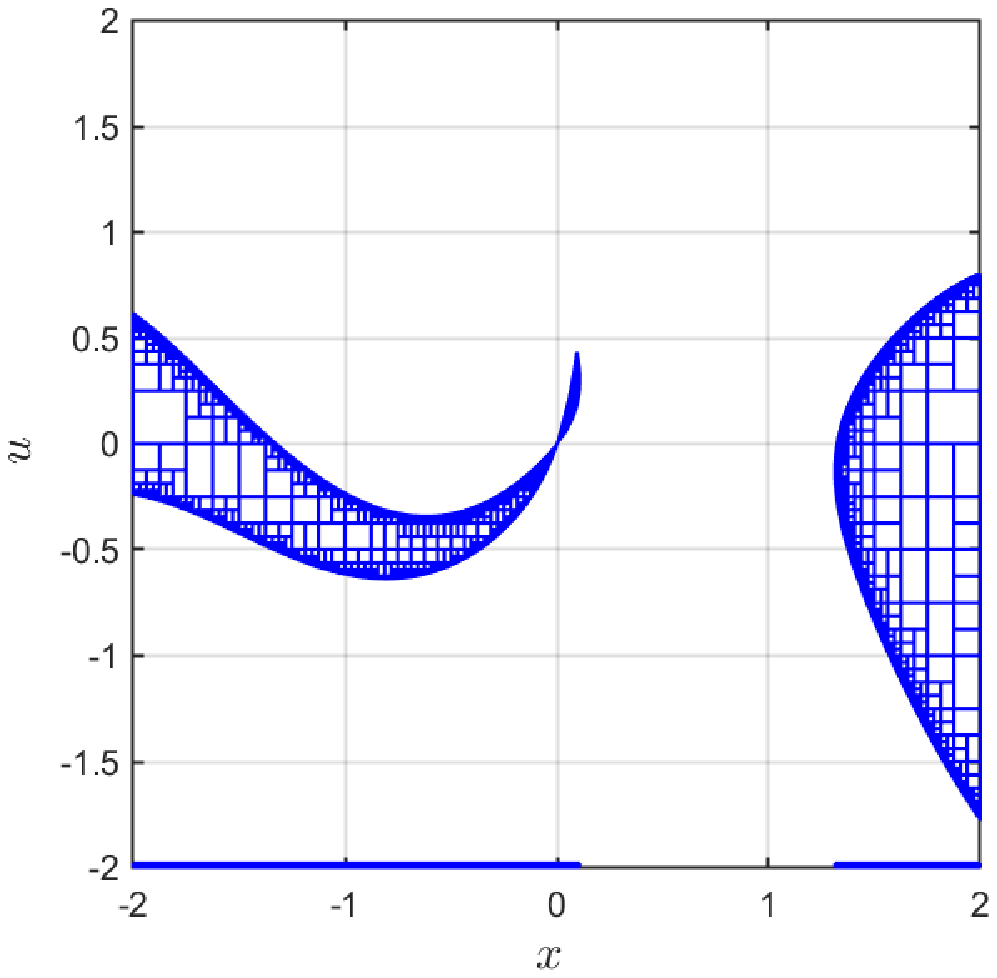}
		\includegraphics[width=0.29\textwidth]{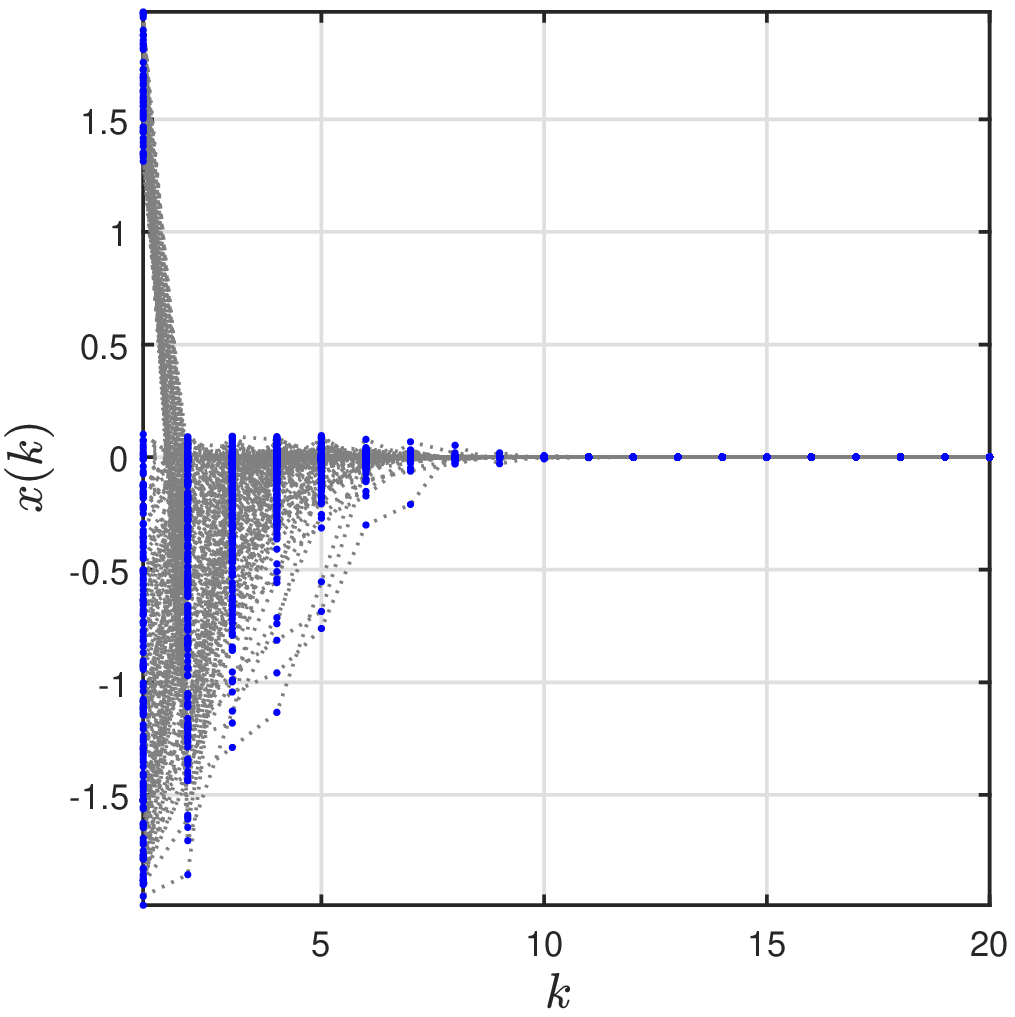}
		\includegraphics[width=0.29\textwidth]{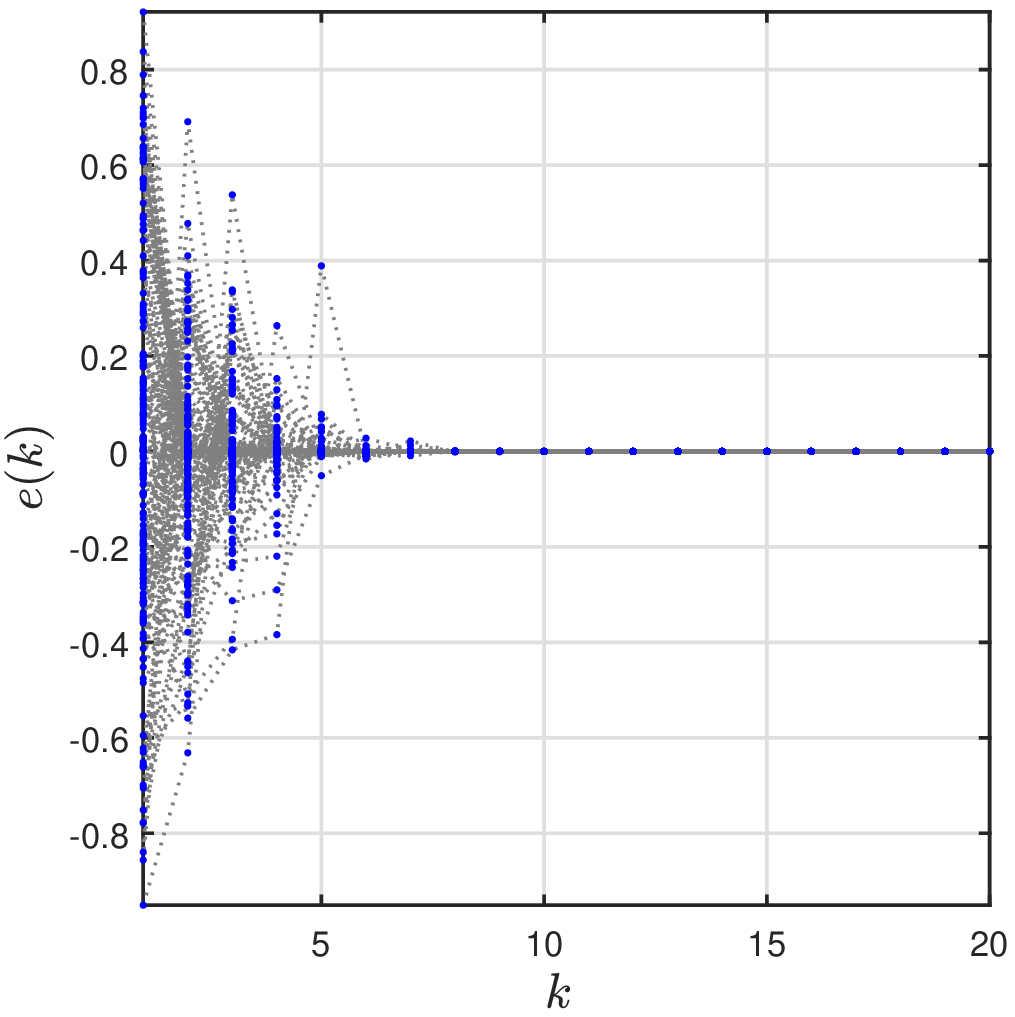}\\
		\parbox[c]{0.29\textwidth}{\footnotesize \centering (a)}
		\parbox[c]{0.29\textwidth}{\footnotesize \centering (b)}
		\parbox[c]{0.29\textwidth}{\footnotesize \centering (c)}
		\includegraphics[width=0.295\textwidth]{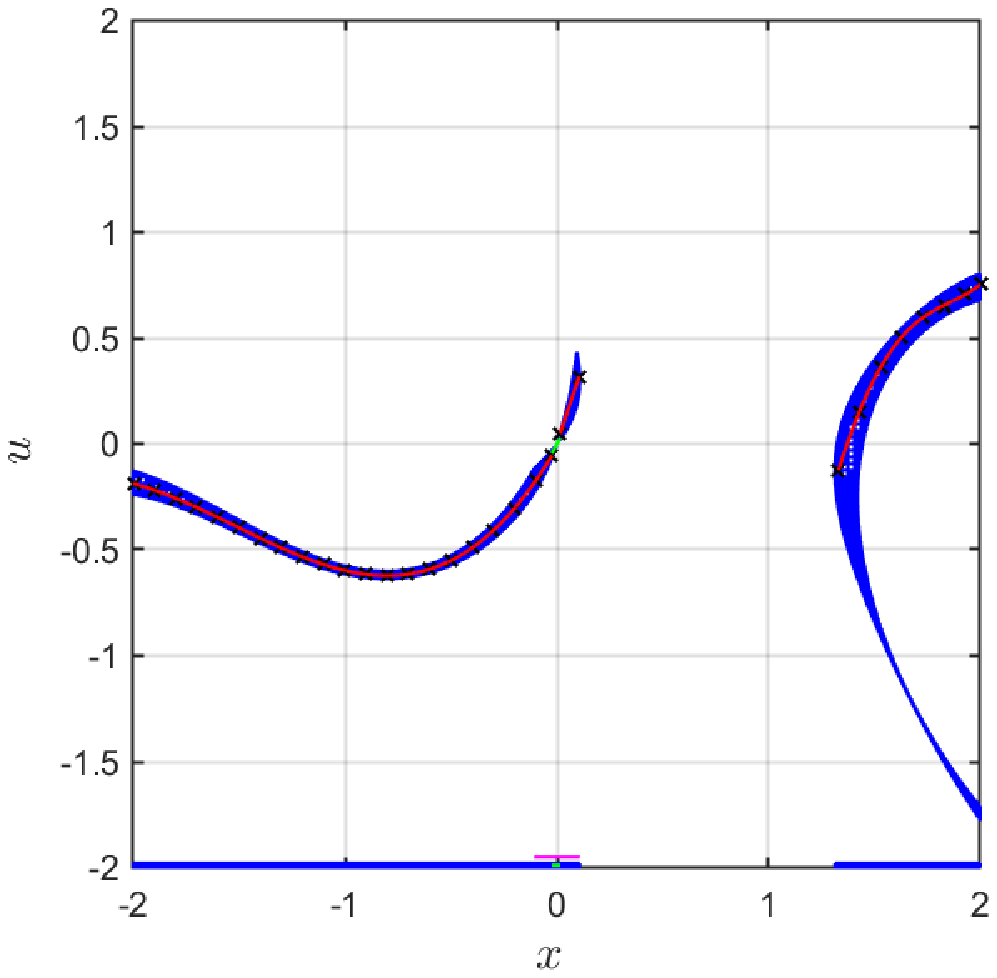}
		\includegraphics[width=0.29\textwidth]{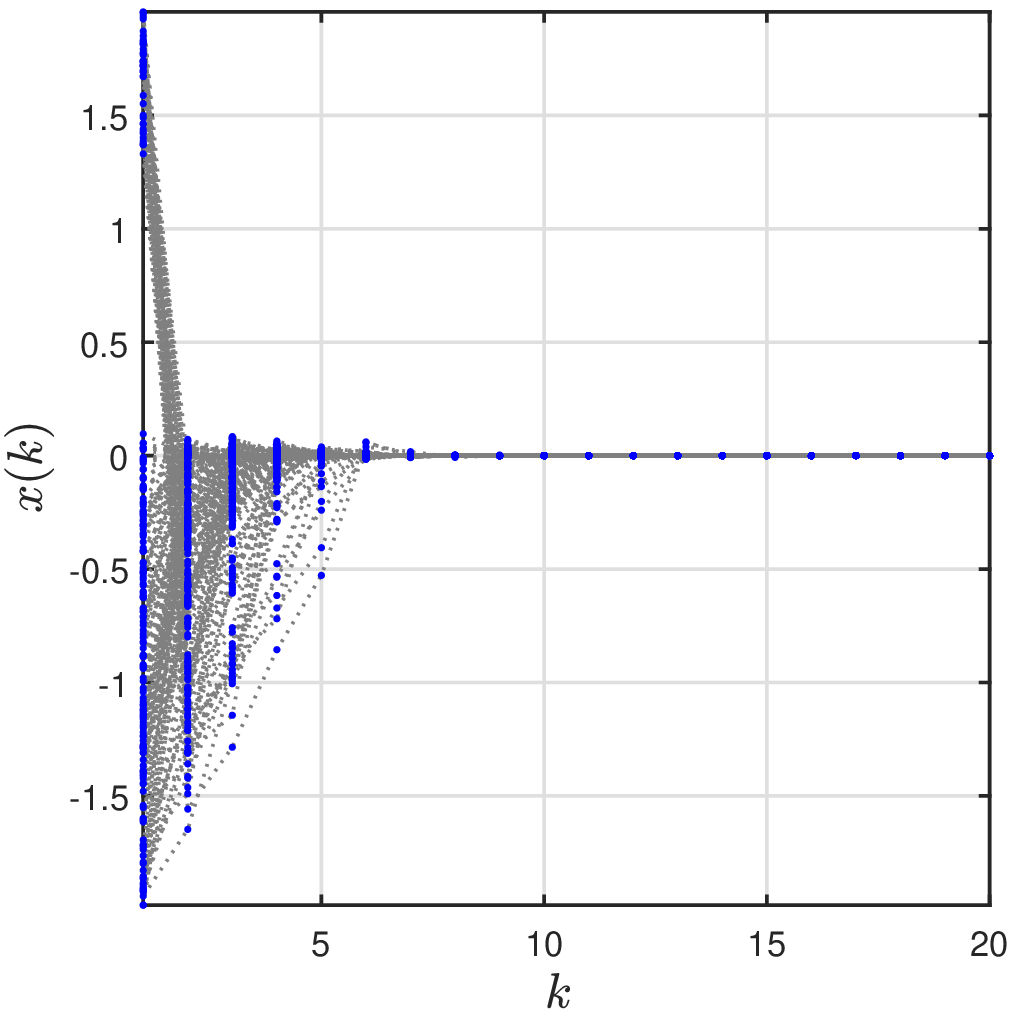}
		\includegraphics[width=0.29\textwidth]{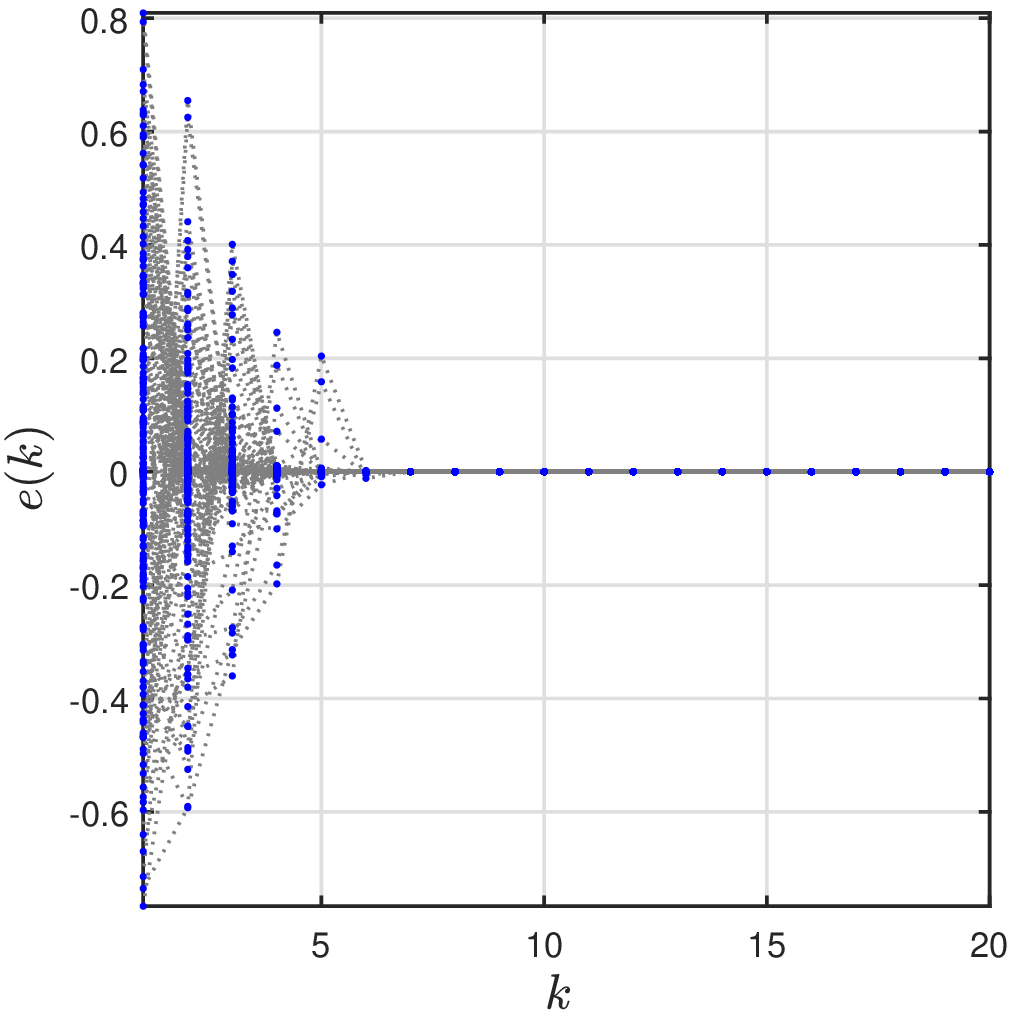}\\
		\parbox[c]{0.29\textwidth}{\footnotesize \centering (d)}
		\parbox[c]{0.29\textwidth}{\footnotesize \centering (e)}
		\parbox[c]{0.29\textwidth}{\footnotesize \centering (f)}
		\caption{(a) Inner approximation $\hat{\mathbb{W}}_{\mathrm{N}}(L)$ of the negative-definite set $\mathbb{W}_{\mathrm{N}}(L)$. (b) State trajectories of closed-loops $x(k)$. (c) Model error trajectories $e(k)$. (d) $\hat{\mathbb{W}}_{\mathrm{N\&I}}(L)$ with controller $\tilde{\mu}$, the estimate of RDOA  $\mathrm{proj}(\hat{\mathbb{W}}_{\mathrm{N}}(L))$. (e) State trajectories of closed-loops $x(k)$ of $\tilde{\mu}$. (f) Model error trajectories $e(k)$ of $\tilde{\mu}$.}
		\label{fig:exmp:x^2}
	\end{center}
\end{figure*}

The given Lyapunov function is selected as $L(x) = x^2$. In Algorithm~\ref{alg:est_W_{NI}_(L)}, the parameter $\epsilon$ is selected as $\epsilon = 0.0001$. The inner approximation $\hat{\mathbb{W}}_{\mathrm{N}}(L)$ of the negative-definite set $\mathbb{W}_{\mathrm{N}}(L)$ is shown in Fig.~\ref{fig:exmp:x^2} (a) denoted by blue boxes. The projection of $\mathbb{W}_{\mathrm{N}}(L)$ along the control space $\mathrm{proj}(\hat{\mathbb{W}}_{\mathrm{N}}(L)) = [-2,-1.25\times10^{-4}]  \cup [-1.25\times10^{-4},0.105] \cup [1.315,2]$ is also shown in  Fig.~\ref{fig:exmp:x^2} (a) denoted by blue line segments in $x$-axis. The inner approximation $\hat{\mathbb{W}}_{\mathrm{N\&I}}(L)$ of the negative-definite and invariant set $\mathbb{W}_{\mathrm{N\&I}}(L)$ is shown in Fig.~\ref{fig:exmp:x^2} (d) denoted by blue boxes. The projection of $\mathbb{W}_{\mathrm{N\&I}}(L)$ along the control space $\mathrm{proj}(\hat{\mathbb{W}}_{\mathrm{N\&I}}(L)) = [-2,-1.25\times10^{-4}] \cup [-1.25\times10^{-4}, 0.105] \cup [1.315,2]$ is also shown in  Fig.~\ref{fig:exmp:x^2} (d) denoted by blue line segments in $x$-axis. The small neighborhood $\mathbb{X}_0$ in \eqref{eq:controller} of the origin is $[-0.03, 0.03]$ as shown in Fig.~\ref{fig:exmp:x^2} (d) denoted by the green line segment in $x$-axis. Hence, the estimate of the closed-loop DOA is $\mathrm{proj}(\hat{\mathbb{W}}_{\mathrm{N\&I}}(L)) \cup \mathbb{X}_0 = [-2,0.105] \cup [1.315,2]$. The result is broader in comparison with the other methods in~\cite{Li:2014_79,Swiatlak:2015} which the invariant of the estimate of the closed-loop RDOA is guaranteed by the level-set of $\mathbb{X}_{\rm ls}(L,0.0117)$, the result is $[-0.108,0.108]$ as shown in Fig.~\ref{fig:exmp:x^2} (d) denoted by the magenta line segment in $x$-axis.

The linear controller in \eqref{eq:controller} is $u = 1.8649x$ denoted by the green straight line through the origin in Fig.~\ref{fig:exmp:x^2} (d). To verify whether all controllers belonging to $\hat{\mathbb{W}}_{\mathrm{N\&I}}(L)$ can stabilize the plant, Fig.~\ref{fig:exmp:x^2} (b) also shows 200 state trajectories of the closed-loop, whose initial states are drawn from the uniform distribution on $[-2,0.105] \cup [1.315,2]$. The trajectories of actual  model error $e(k)\in[- \delta(x(k),u(k)),+\delta(x(k),u(k))]$ are shown in Fig.~\ref{fig:exmp:x^2} (c). Here, the linear controller is $u = 1.8649x$, while the output of the nonlinear controller is drawn from the uniform distribution on $\mathbb{U}(x) \subset \mathbb{R}$. For the given $x \in \mathbb{R}$, $\mathbb{U}(x)$ is defined as $\mathbb{U}(x) = \big\{u \in \mathbb{R} | (x;u) \in \hat{\mathbb{W}}_{\mathrm{N\&I}}(L) \big\}$.

In order to find the nonlinear controller $\tilde{\mu}$ in \eqref{eq:controller}, we select a training data set denoted by black $'\times '$s in Fig.~\ref{fig:exmp:x^2} (d). Then, $\tilde{\mu}$ is obtained using Gaussian processes regression, denoted by red line in Fig.~\ref{fig:exmp:x^2} (d). Fig.~\ref{fig:exmp:x^2} (e) shows 200 state trajectories of the closed-loop, whose initial states are drawn from the uniform distribution on $[-2,0.105] \cup [1.33,2]$.  The trajectories of actual  model error $e(k)\in[- \delta(x(k),u(k)),+\delta(x(k),u(k))]$ are shown in Fig.~\ref{fig:exmp:x^2} (f). We see that both state  and model error trajectories of two different controllers  converge to the origin.

State trajectories which are shown in Fig.~\ref{fig:exmp:x^2} (b,e) also verify that $\hat{\mathbb{W}}_{\mathrm{N\&I}}(L) \cup \{(x;u)|u = 1.8649x, x \in \mathbb{X}_0 \}$ is invariant for the plant. We know that $\mathrm{proj}(\hat{\mathbb{W}}_{\mathrm{N\&I}}(L) \cup \{(x;u)|u = 1.8649x, x \in \mathbb{X}_0 \}) = \mathrm{proj}(\hat{\mathbb{W}}_{\mathrm{N\&I}}(L)) \cup \mathbb{X}_0 = [-2,0.105] \cup [1.315,2]$. From the figures, we see that there is no state in $(0.105,1.315) \subset \mathbb{R}$.

\subsection{Stabilization with closed-loop DOA enlargement}
The positive-definite function set in optimization problem \eqref{eq:optim} is selected as
\begin{equation}
\mathfrak{L}_{1,4} = = \Big\{L \in \mathfrak{R}_{1,4} \Big| L(x) = (x;x^2)^TP^TP(x;x^2), x \in \mathbb{R}\Big\} \nonumber
\end{equation}
with the parameters $P \in \mathbb{R}^{2 \times 2}$. The optimization problem \eqref{eq:optim_P} is solved through the particle swarm optimization method and the solution $L^\ast(x) = 1.5327x^4 +2.3121x^3 + 1.1286x^2$ is obtained. 

\begin{figure*}
	\begin{center}
		\includegraphics[width=0.295\textwidth]{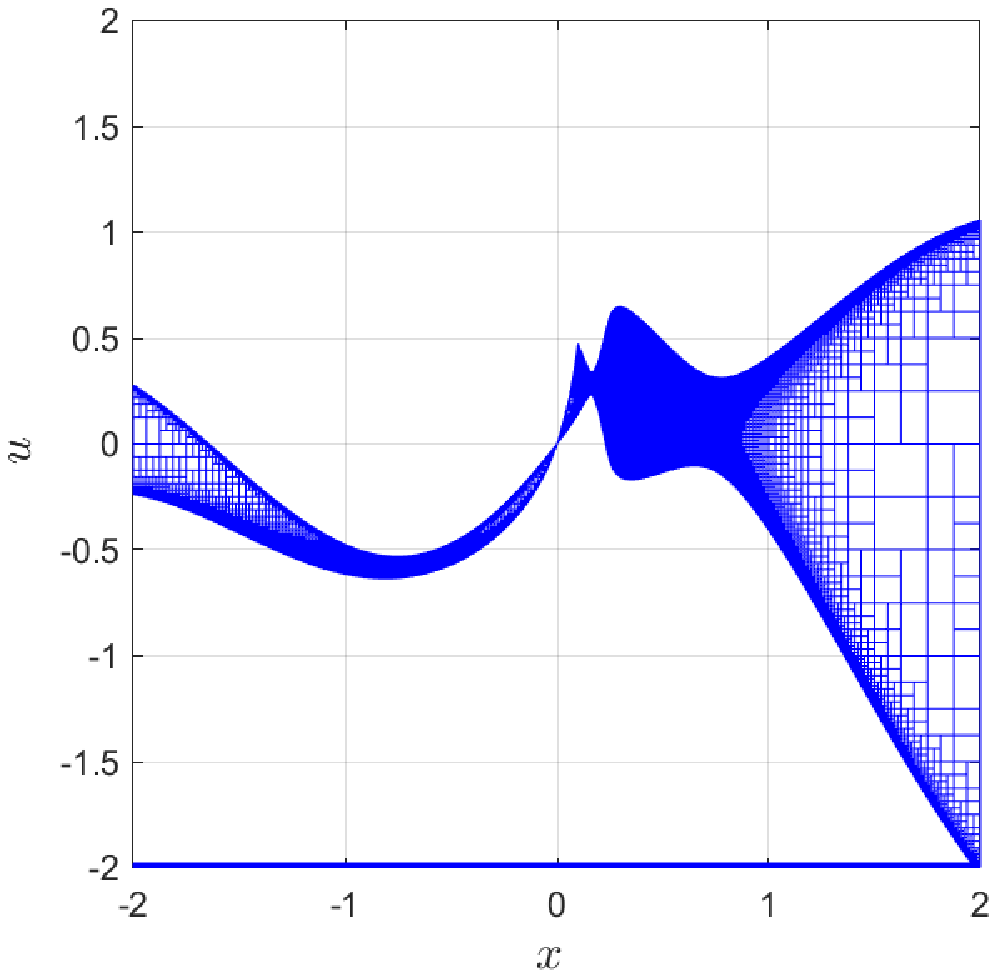}
		\includegraphics[width=0.29\textwidth]{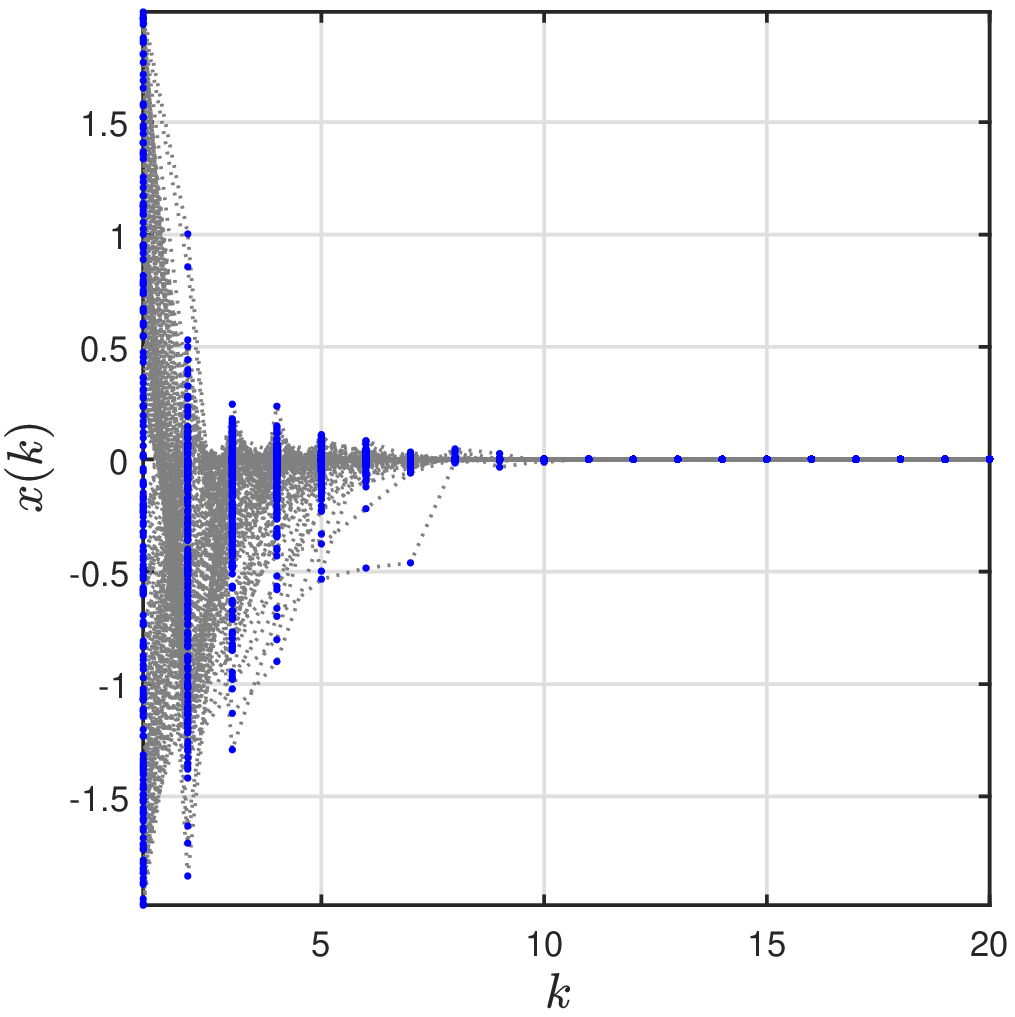}
		\includegraphics[width=0.29\textwidth]{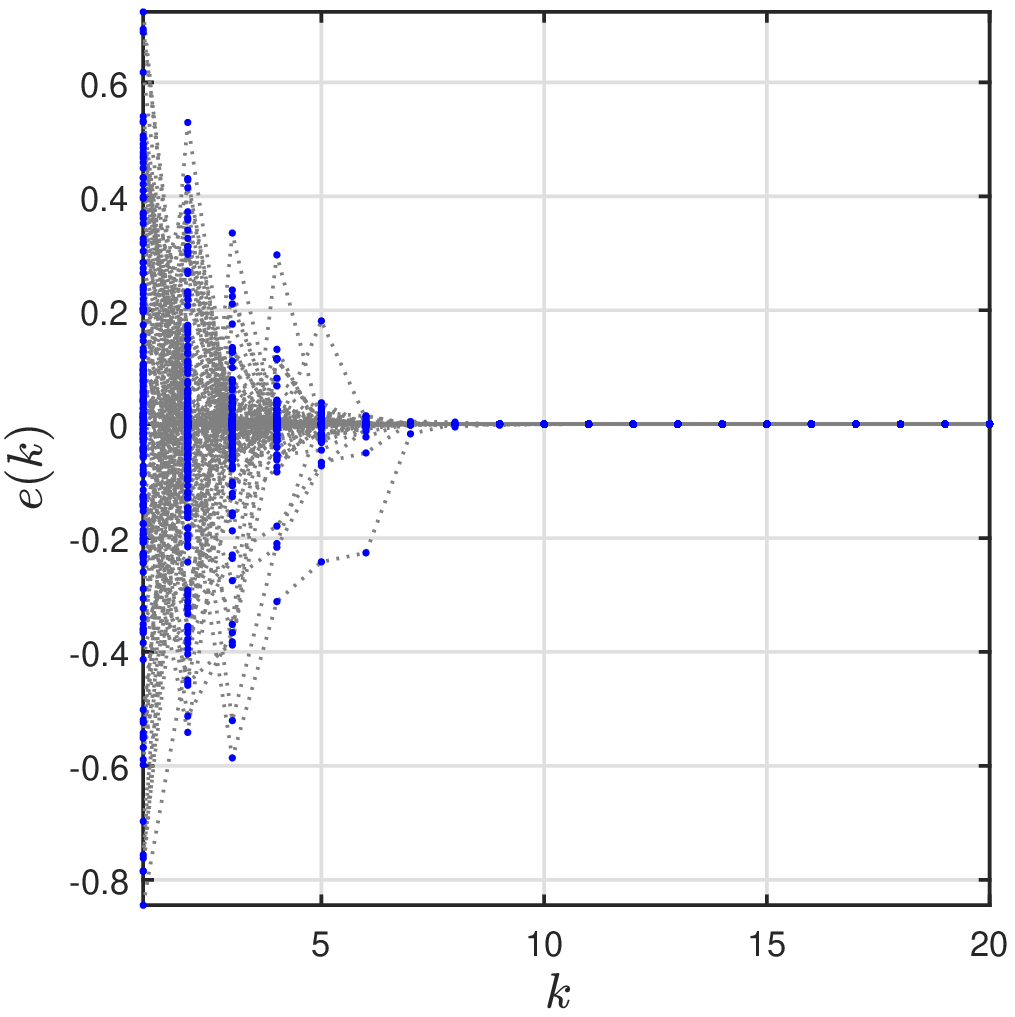}\\
		\parbox[c]{0.29\textwidth}{\footnotesize \centering (a)}
		\parbox[c]{0.29\textwidth}{\footnotesize \centering (b)}
		\parbox[c]{0.29\textwidth}{\footnotesize \centering (c)}\\
		\includegraphics[width=0.295\textwidth]{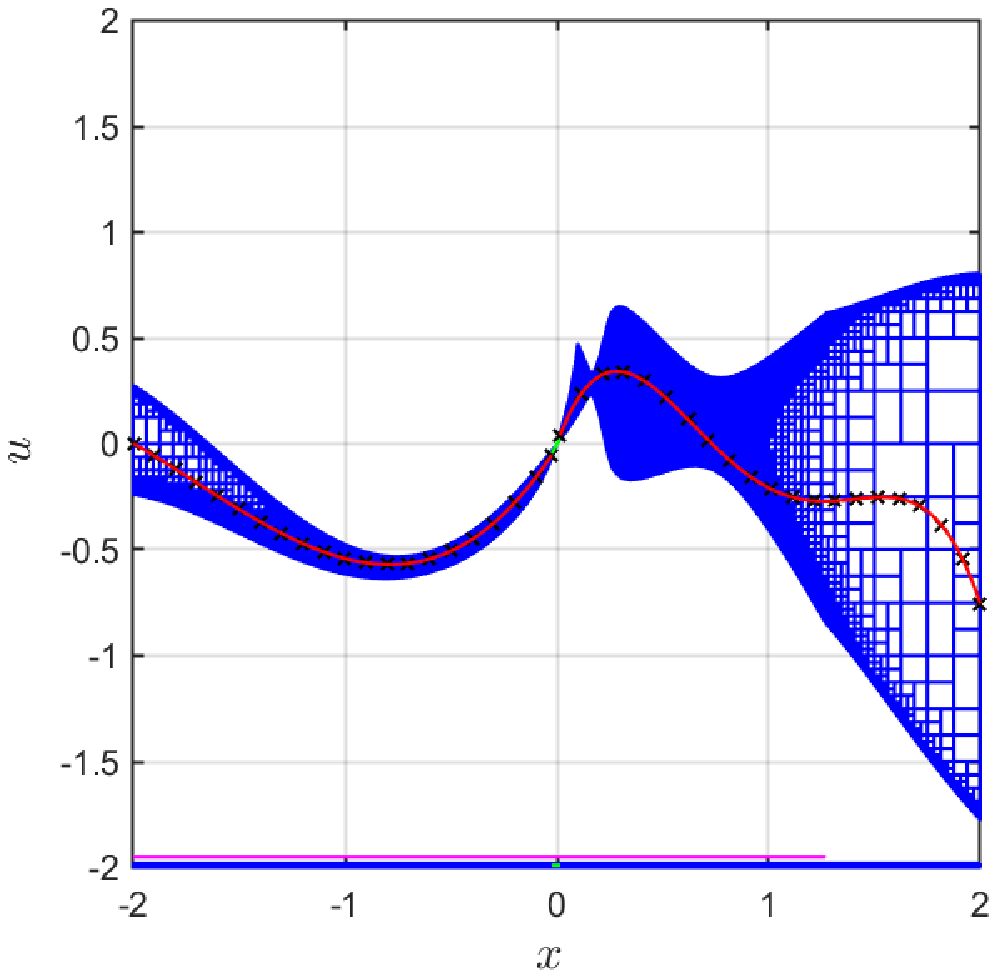}
		\includegraphics[width=0.29\textwidth]{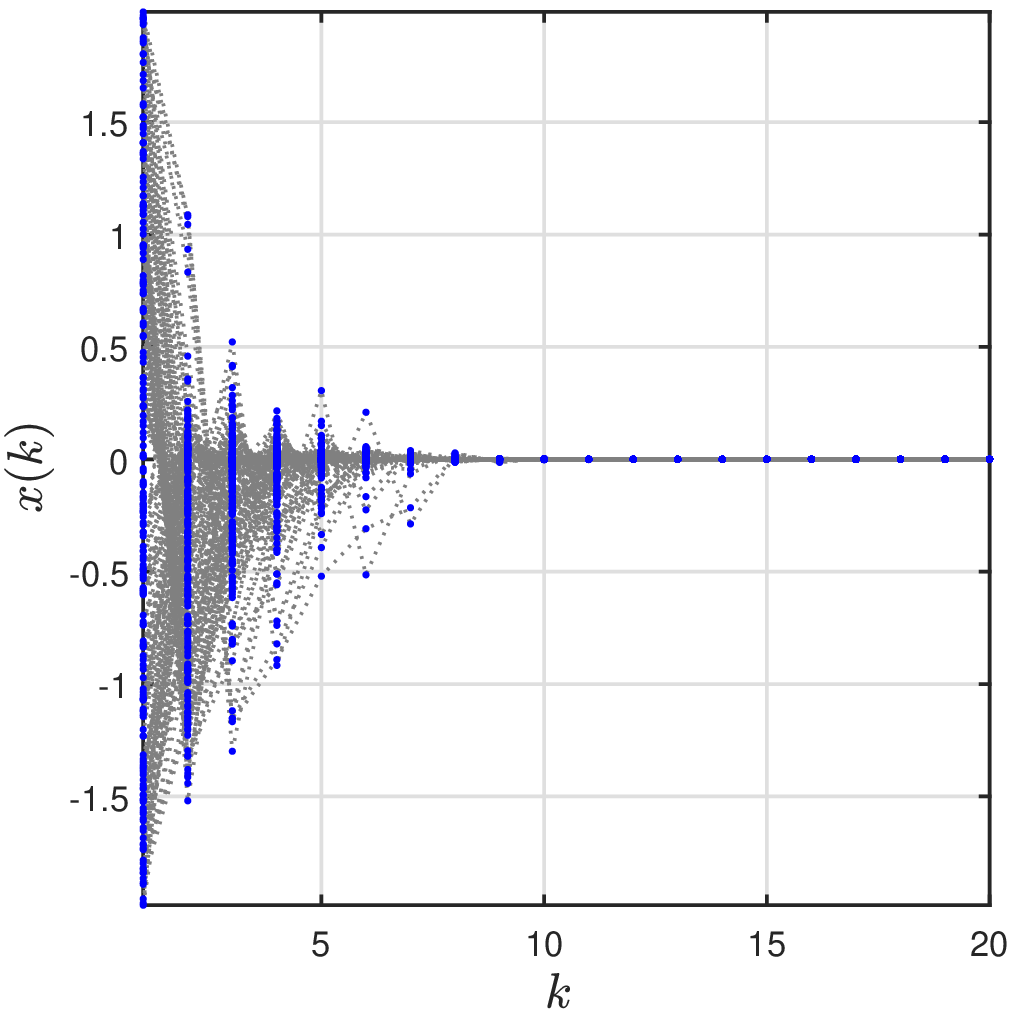}
		\includegraphics[width=0.29\textwidth]{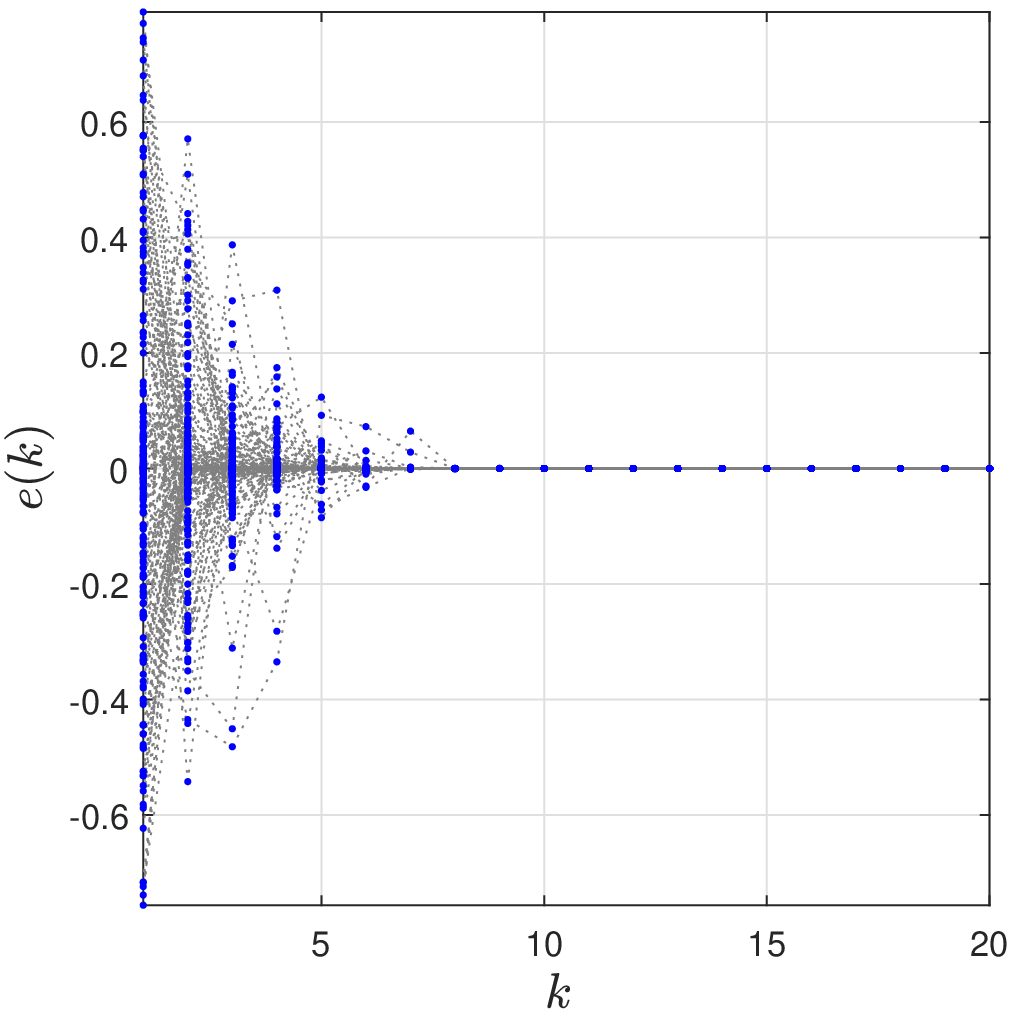}\\
		\parbox[c]{0.29\textwidth}{\footnotesize \centering (d)}
		\parbox[c]{0.29\textwidth}{\footnotesize \centering (e)}
		\parbox[c]{0.29\textwidth}{\footnotesize \centering (f)}
		\caption{(a) Inner approximation $\hat{\mathbb{W}}_{\mathrm{N}}(L^\ast)$ of the negative-definite set $\mathbb{W}_{\mathrm{N}}(L^\ast)$. (b) State trajectories of closed-loops $x(k)$. (c) Model error trajectories $e(k)$. (d) $\hat{\mathbb{W}}_{\mathrm{N\&I}}(L^\ast)$ with controller $\tilde{\mu}$, the estimate of RDOA  $\mathrm{proj}(\hat{\mathbb{W}}_{\mathrm{N}}(L^\ast))$. (e) State trajectories of closed-loops $x(k)$ of $\tilde{\mu}$. (f) Model error trajectories $e(k)$ of $\tilde{\mu}$.}
		\label{fig:exmp:opt}
	\end{center}
\end{figure*}

The inner approximation $\hat{\mathbb{W}}_{\mathrm{N}}(L^\ast)$ of the negative-definite set $\mathbb{W}_{\mathrm{N}}(L^\ast)$ is shown in Fig.~\ref{fig:exmp:opt} (a) denoted by blue boxes. The inner approximation $\hat{\mathbb{W}}_{\mathrm{N\&I}}(L^\ast)$ of the negative-definite and invariant set $\mathbb{W}_{\mathrm{N\&I}}(L^\ast)$ is shown in Fig.~\ref{fig:exmp:opt} (d) denoted by blue boxes. The projection of $\mathbb{W}_{\mathrm{N\&I}}(L^\ast)$ along the control space $\mathrm{proj}(\hat{\mathbb{W}}_{\mathrm{N\&I}}(L^\ast)) = [-2,-1.25\times10^{-4}] \cup [1.25\times10^{-4},2]$ is also shown in  Fig.~\ref{fig:exmp:opt} (d) denoted by blue line segments in $x$-axis. The small neighborhood $\mathbb{X}_0$ in \eqref{eq:controller} of the origin is $[-0.03,0.03]$ as shown in Fig.~\ref{fig:exmp:opt} (d) denoted by the green line segment in $x$-axis. Hence, the estimate of the closed-loop DOA is $\mathrm{proj}(\hat{\mathbb{W}}_{\mathrm{N\&I}}(L^\ast)) \cup \mathbb{X}_0 = [-2,2]$. The result is broader in comparison with another method in~\cite{Li:2014_79} which the invariant of the estimate of the closed-loop RDOA is guaranteed by the level-set of $\mathbb{X}_{\rm ls}(L,10.5408)$, the result is $[-2,1.27]$ as shown in Fig.~\ref{fig:exmp:opt} (d) denoted by the magenta line segment in $x$-axis.

The linear controller in \eqref{eq:controller} is $u = 1.8649x$ denoted by the green straight line through the origin in Fig.~\ref{fig:exmp:opt} (d), while the output of the nonlinear controller is drawn from the uniform distribution on $\mathbb{U}(x) \subset \mathbb{R}$. To verify whether all controllers belonging to $\hat{\mathbb{W}}_{\mathrm{N\&I}}(L^\ast)$ can stabilize the plant, Fig.~\ref{fig:exmp:opt} (b) also shows 200 state trajectories of the closed-loop, whose initial states are drawn from the uniform distribution on $[-2,2]$. The trajectories of actual  model error $e(k)$ are shown in Fig.~\ref{fig:exmp:opt} (c).

For the given $x \in \mathbb{R}$, $\mathbb{U}(x)$ is defined as $\mathbb{U}(x) = \big\{u \in \mathbb{R} | (x;u) \in \hat{\mathbb{W}}_{\mathrm{N\&I}}(L^\ast) \big\}$.In order to find the nonlinear controller $\tilde{\mu}$ in \eqref{eq:controller}, we select a training data set denoted by black $'\times '$s in Fig.~\ref{fig:exmp:opt} (d). Then, $\tilde{\mu}$ is obtained using Gaussian processes regression, denoted by red line in Fig.~\ref{fig:exmp:opt} (d). Fig.~\ref{fig:exmp:opt} (e) shows 200 state trajectories of the closed-loop, whose initial states are drawn from the uniform distribution on $[-2,2]$. The trajectories of actual  model error $e(k)\in[- \delta(x(k),u(k)),+\delta(x(k),u(k))]$ are shown Fig.~\ref{fig:exmp:opt} (f). We see that all state trajectories converge to the origin.

\section{Conclusion}
We have proposed an interval analysis approach to find an unstructured robust controller set for the nonlinear discrete-time systems with nonlinear uncertainties. The results presented above shown the interest of the proposed approach: 1)For a given Lyapunov function, the estimate of RDOA is successfully obtained  with respect to the RNIS-SC and is boarder than the Lyapunov function under the level set restriction. 2)For the different Lyapunov functions, the RDOA is totally different. After the solvable optimization problem is formulated, the RDOA of RNIS-SC is enlarged by selecting an appropriate Lyapunov function from a positive-definite function set. According to the results, the enlarged RDOAs are much better and verify the effectiveness of the proposed method. 

For discrete-time systems with fixed modeling error bound at every time instant K, the next state $x(k+1)$ is always in a fixed interval determined by the true plant and uncertainty. Therefore, our presented RNIDEVIA algorithm based on interval arithmetic is very competent to address this issue. Moreover, the estimate of RDOA in interval form can guarantee the rigor in the computed result, that is why we name our scheme as an interval-driven method.

\appendix
\section{Computation Detailes of Simulation} \label{Appendix}

The considered example considers a scalar state and a scalar control input. The pseudocode for estimating the robust negative-definite set $\hat{\mathbb{W}}_{\rm N}$ for the given plant set is shown in Pesudocode~\ref{alg:RNDDA}. The whole code is present in $\rm https://github.com/CharlieLuuke/IDRCS$. Computation time results are shown in Tabel~\ref{Tabel:nega}.
~\\

\setcounter{algorithm}{0}
\begin{balgorithm} \label{alg:RNDDA}
	\caption{The procedure of obtaining $\hat{\mathbb{W}}_{\rm N}$ } 
	\textbf{Inputs}:
	
	\ \ -\ Lyapunov-definite function $L: \mathbb{R}^n \to \mathbb{R}$;
	
	\ \ -\ Nominal model $f(x,u)$.
	
	\ \ -\ Error bound $\delta(x,u)$.
	
	\ \ -\ Interested region $\mathbb{W}_0$.		
	
	\textbf{Outputs}:
	
	\ \ -\ Robust negative definite set $\mathbb{W}_R$.
	
	\textbf{Parameters}:
	
	\ \ -\ Desired accuracy $\epsilon=1e^{-4}$; Coefficient of Lyapunov function $a,b,c$.
	
	\ \ -\ $\mathbb{W}_{list}=\emptyset,\mathbb{W}_{in}=\emptyset,\mathbb{W}_{bou}=\emptyset,\mathbb{W}_{out}=\emptyset$.
	
	\textbf{Steps}: 
	
	\begin{algorithmic}[1]
		
		\State Construct function $V=a(f(x,u)\pm\delta(x,u))^2 + b(f(x,u)\pm\delta(x,u))^3+c(f(x,u)\pm\delta(x,u))^4-(ax^2+bx^3+cx^4)$
		\State $\mathbb{W}_{\mathrm{list}} :=[\mathbb{W}_{\mathrm{list}},\mathbb{W}_0]$
		\While {$\mathbb{W}_\mathrm{list} \neq \emptyset$}
		\State $\mathbb{W}_{do}:={\mathbb{W}}_{\mathrm{list}}[end]$
		\State $\bar{\mathbb{X}}_{f+\delta}=V_{(f+\delta)}.eval(\mathbb{W}_\mathrm{do})$
		\State $\bar{\mathbb{X}}_{f-\delta}=V_{(f-\delta)}.eval(\mathbb{W}_\mathrm{do})$
		\State Remove ${\mathbb{W}}_{\mathrm{list}}[end]$
		\If {$max(\bar{\mathbb{X}}_{f+\delta})<0$ and $ max(\bar{\mathbb{X}}_{f-\delta})<0$}
		\State $\mathbb{W}_{\mathrm{in}} :=[\mathbb{W}_{\mathrm{in}},\mathbb{W}_\mathrm{do}]$
		\ElsIf {$min(\bar{\mathbb{X}}_{f+\delta})>0$ or $ min(\bar{\mathbb{X}}_{f-\delta})>0$}
		\State $\mathbb{W}_{\mathrm{out}} :=[\mathbb{W}_{\mathrm{out}},\mathbb{W}_\mathrm{do}]$
		\ElsIf {$d(\mathbb{W}_\mathrm{do}) < \epsilon$}
		\State $\mathbb{W}_{\mathrm{bou}} :=[\mathbb{W}_{\mathrm{bou}},\mathbb{W}_\mathrm{do}]$
		\Else
		\State Bisect box $\mathbb{W}_\mathrm{do}$
		\State Add the new boxes to set $\mathbb{W}_{\mathrm{list}}$
		\EndIf
		\EndWhile
		\State \textbf{return} ${\mathbb{W}}_{\mathrm{in}}$%, Z_{\mathrm{bou}}, Z_{\mathrm{out}}$
		
	\end{algorithmic}
\end{balgorithm}

\begin{table}[h]   
	\tiny  
	\caption{Comparison of Run Times in Two Platforms}  
	\begin{tabular}{|l|c|c|} 
		\hline
		\diagbox {L. F.}{Time (s)}{Platform} & i5 3.2 GHZ & Athlon 2.7 GHZ \\  
		\hline  
		$x^2$ & 21.16 & 43.69 \\  
		\hline  
		$L^\ast(x)$ & 187.23 & 346.72 \\  
		\hline  
		
	\end{tabular} 
    \label{Tabel:nega}   
\end{table}

\begin{ack}                               
	The authors would like to thank Dr. Yinan Li for helpful discussions. All computations were performed with the Pyibex~\cite{Desrochers:DB}. 
	
\end{ack}

\bibliographystyle{plain}        % Include this if you use bibtex
\bibliography{stab_interval}       % and a bib file to produce the
% bibliography (preferred). The
% correct style is generated by
% Elsevier at the time of printing.

%\section{Proof of *} \label{app:proof:thm:robust_stab}
%
%\begin{pf}
%	***
%\end{pf}

%\section{Some Latin vocabulary}         % Sections and subsections are supported
                                        % in the appendices.
                
\end{document}